\documentclass[12pt,a4paper]{article}
\usepackage{amssymb,amsthm,amsmath}
\usepackage[english]{babel}
\usepackage{fullpage}

\title{The Isospectral Dirac Operator on the $4$-dimensional \\ [6pt] Orthogonal Quantum  Sphere}

\date{13th February 2008}

\author{~\\
\large{Francesco D'Andrea$^1$, Ludwik D\c{a}browski$^1$, 
 	Giovanni Landi$^2$} \\ [10pt]
 \normalsize{$^1$ Scuola Internazionale Superiore di Studi Avanzati,}\\
 \normalsize{Via Beirut 2-4, I-34014, Trieste, Italy} \\ [10pt]
 \normalsize{$^2$ Dipartimento di Matematica e Informatica, Universit{\`a} di Trieste,} \\
\normalsize{Via  Valerio 12/1, I-34127, Trieste, Italy} \\
\normalsize{and INFN, Sezione di Trieste, Trieste, Italy}
\\~\\
}

\allowdisplaybreaks[4]
\newenvironment{prova}{\begin{proof}\parskip=0mm\parindent=0in}{\end{proof}}
\numberwithin{equation}{section}

\newtheorem{thm}{Theorem}[section]
\newtheorem{lemma}[thm]{Lemma}
\newtheorem{prop}[thm]{Proposition}
\newtheorem{df}[thm]{Definition}

\newcommand{\A}{\mathcal{A}}
\newcommand{\B}{\mathcal{B}}
\newcommand{\U}{\mathcal{U}}
\newcommand{\D}{D\mkern-11.5mu/\,}
\newcommand{\HH}{\mathcal{H}}
\newcommand{\mU}{\mathcal{U}}
\newcommand{\Uq}{U_q(so(5))}
\newcommand{\C}{\mathbb{C}}
\newcommand{\R}{\mathbb{R}}
\newcommand{\N}{\mathbb{N}}
\newcommand{\Z}{\mathbb{Z}}

\newcommand{\sqbn}[2]{\Big[\textrm{\footnotesize$\!\!\begin{array}{c}#1 \\ #2\end{array}\!\!$}\Big]}
\newcommand{\ket}[1]{\left|#1\right>}
\newcommand{\inner}[1]{\left<#1\right>}
\newcommand{\ma}[1]{\textrm{\scriptsize$\left(\!\begin{array}{cccc}#1\end{array}\!\right)$}}

\newcommand{\az}{\triangleright}
\newcommand{\tr}{\mathrm{Tr}}
\newcommand{\de}{\mathrm{d}}
\newcommand{\op}{\mathrm{OP}^{-\infty}}
\newcommand{\opz}{\mathrm{OP}^0}
\newcommand{\nint}{\int\mkern-19mu-\;}
\newcommand{\openone}{\leavevmode\hbox{\small1\kern-3.8pt\normalsize 1}}
\newcommand{\kkett}[1]{\,\left|\mkern-1mu\left|\smash[b]{\smash[t]{#1}}\right>\mkern-3mu\right>}
\newcommand{\oh}{\smash[t]{\smash[b]{\tfrac{1}{2}}}}
\newcommand{\idJ}{\mathcal{I}}

\begin{document}

\maketitle

\bigskip

\thispagestyle{empty}
                                                                               
\begin{abstract}
\noindent 
Equivariance under the action of $\Uq$ is used to compute the left
regular and (chiral) spinorial representations of the algebra of
the orthogonal quantum  $4$-sphere $S^4_q$. These representations are
the constituents of a spectral triple on $S^4_q$ with a Dirac operator
which is isospectral to the canonical one on the
round sphere $S^4$ and which then gives $4^+$-summability. Non-triviality of
the geometry is proved by pairing the associated Fredholm module with an
`instanton' projection.
We also introduce a real structure which satisfies all required properties
modulo smoothing operators.
\end{abstract}

\vfill


\pagebreak

\section{Introduction}

The recent constructions of spectral triples -- with the consequent
analysis of the corresponding spectral geometry -- for the manifold
of the quantum $SU(2)$ group in~\cite{CP,ConSU,DLSvSV,vSDLSV}
and for its quantum homogeneous spaces (the Podle\'s spheres)
in~\cite{DS,DLPS,DD,DDLW}, have provided a number of examples
showing that a marriage between noncommutative geometry and quantum
groups theory is indeed possible. A common feature of most of these
examples is that the dimension spectrum is the same as in the commutative
($q=1$) limit. Furthermore, with the only known exception of the
$0^+$-summable `exponential' spectral triple on the standard Podle\'s sphere
given in~\cite{DS}, in order to have a real spectral triple one is forced 
to weaken the usual requirements that the real structure should satisfy.

It is then only natural to try and construct additional explicit examples wondering
in particular if these properties are common to all quantum spaces or are
rather coincidences which happen for low dimensional examples (all related
to the quantum group $SU_q(2)$).
In this paper we present an example in `dimension four' given by a spectral
triple on the orthogonal quantum  sphere $S^4_q$ which is isospectral to the
canonical spectral triple on the classical sphere with the round metric.
There exists also a real structure which satisfies  all required properties
modulo an ideal of smoothing operators.

There are a few reasons why in dimension greater than or equal to four 
the orthogonal quantum  sphere $S^4_q$ is most interesting to study.
Firstly, all the relevant irreducible
representations of the symmetry algebra $\Uq$ are known~\cite{Cha} and 
both the algebra $\A(S^4_q)$ of polynomial functions as well as the
modules of chiral spinors carry representations of $\Uq$
which are multiplicity free. Secondly, the spectrum of the Dirac operator $\D$
for the round metric on the undeformed sphere $S^4$ is known~\cite{Trau,CH}.
All this allows us to apply the already tested methods of isospectral
deformations and 
to construct an $\Uq$-equivariant spectral triple on $S^4_q$.

The sphere $S^4_q$ could also be relevant for noncommutative physical models.
In particular, on  $S^4_q$ there is a canonical `instantonic vector bundle'~\cite{LH}
and the study of the noncommutative geometry of $S^4_q$ could be a first step for
the construction of $SU_q(2)$ instantons on this space.

\bigskip

In Sect.~\ref{se:pre} we recall all generalities about spectral triples that we need.
We give also some properties of finitely generated projective modules over algebras
having quantum group symmetries. The rest of the paper is organized as follows.
Sections \ref{sec:rep} and \ref{sec:S4q} are devoted to the symmetry 
Hopf algebra $\Uq$ and its fundamental $*$-algebra module,
the orthogonal quantum  sphere $S^4_q$. In Sect.~\ref{sec:cs}
we describe the $\A(S^4_q)$-modules of chiral spinors over $S^4_q$.
Sect.~\ref{sec:eqreps} is devoted to the
left regular representation of the algebra $\A(S^4_q)$ of polynomial functions over
$S^4_q$ and to the representations of $\A(S^4_q)$ which in the $q=1$ limit correspond
to the modules of chiral spinors. These representations are $\Uq$-equivariant, that is
they correspond to representations of the crossed product algebra $\A(S^4_q)\rtimes\Uq$.
In Section \ref{sec:D} we use the isospectral Dirac operator to construct a spectral
triple on $S^4_q$; it will be $\Uq$-equivariant,  regular, even and of metric
dimension $4$. We also prove that it is non-trivial by pairing the Fredholm module
canonically associated to the spectral triple to an `instanton' projection $e$.
It turns out that the projection $e$ has charge $1$, as in the classical case.
In Sect.~\ref{sec:Sigma}, we compute the part of the dimension spectrum contained
in the right half plane $\{s\in\C\,|\,\mathrm{Re}\,s>2\}$, as well as the top
residue (which in the commutative case is proportional to the integral).
This is done by quotienting by a suitable ideal of `infinitesimals' $\idJ$,
which is larger than smoothing operators.
At the moment we are unable to comment on the part of the dimension spectrum which
is in the left half plane $\mathrm{Re}\,s\leq 2$.
Finally, in Sect.~\ref{sec:real} we produce an equivariant real structure for
which both the `commutant  property' and the `first order condition' are satisfied
modulo the ideal of smoothing operators; this is consonant with the cases of the
manifold of $SU_q(2)$ in \cite{DLSvSV} and of Podle{\'s} spheres in \cite{DLPS,DDLW}.
In fact, we also show that these conditions are much easier to handle modulo the ideal $\idJ$.
 
\section{Some useful preliminaries}\label{se:pre}

In this section, we collect some basic notions concerning equivariant spectral triples.
We also give some general properties of finitely generated projective modules over
algebras having quantum group symmetries. 

\subsection{Generalities about Spectral Triples}\label{app:A}
We start with the notion of finite summable spectral triples~\cite{CBook}.
\begin{df}
A \emph{spectral triple} $(\A,\HH,D)$ is the datum of a complex associative unital
$*$-algebra $\A$, a $*$-representation $\pi:\A\to\B(\HH)$ by bounded operators on a
(separable) Hilbert space $\HH$ and a self-adjoint (unbounded) operator $D=D^*$ such that,
\begin{itemize}
\item $\;\;(D+i)^{-1}$ is a compact operator;
\item $\;\;[D,\pi(a)]$ is a bounded operator for all $a\in\A\,$.
\end{itemize}
\end{df}

\noindent We refer to $D$ as the `generalized' Dirac operator, or the Dirac operator `tout court' and for simplicity we assume that it is invertible.
Usually, the representation symbol $\pi$ is removed when no risk of confusion arises. 

With $n\in\R^+$, $D$ is called \emph{$n^+$-summable} if the operator 
$(D^2 +1)^{-1/2}$ is in the Dixmier ideal ${\cal L}^{n+}(\HH)$. 
We shall also call $n$  the \emph{metric dimension} of the spectral triple.

A spectral triple is called \emph{even} if there exists a grading
$\gamma$, i.e.~a bounded operator satisfying $\gamma=\gamma^*$ and $\gamma^2=1$,
such that the Dirac operator is odd and the algebra is even:
$$
\gamma D+D\gamma=0\;,\qquad a\gamma=\gamma a \;, \quad \forall\;a\in\A\;.
$$

We recall from~\cite{CM} a few analytic properties of spectral triples. 
To the unbounded operator $D$ on $\HH$ one associates an unbounded
derivation $\delta$ on $\B(\HH)$ by,
$$
\delta(a)=[|D|,a]\;,
$$
for all $a\in\B(\HH)$. A spectral triple is called \emph{regular} if
the following inclusion holds,
$$
\A\cup[D,\A]\subset\bigcap\nolimits_{j\in\N}\mathrm{dom}\,\delta^j\;,
$$
and we refer to $\opz:=\bigcap\nolimits_{j\in\N}\mathrm{dom}\,\delta^j$ as
the `smooth domain' of the operator $\delta$. For a regular spectral triple,
the class $\Psi^0$ of pseudodifferential operators of order less than or equal
to zero is defined as the algebra generated by $\bigcup_{k\in\N}\delta^k(\A\cup[D,\A])$.
If the triple has finite metric dimension $n$, the `zeta-type' function 
$$
\zeta_a(s):=\tr_{\HH}(a|D|^{-s})
$$
associated to $a\in\Psi^0$ is defined (and holomorphic) for $s\in\C$ with $\mathrm{Re}\,s>n$
and the following definition makes sense.
\begin{df}
A spectral triple has \emph{dimension spectrum} $\Sigma$ if{}f $\Sigma\subset\C$
is a countable set, for all $a\in\Psi^0$ the function $\zeta_a(s)$ extends to a
meromorphic function on $\C$ with poles as unique singularities, and the union of such
singularities is the set $\Sigma$.
\end{df}

\noindent
If $\Sigma$ is made only of simple poles, the Wodzicki-type residue functional
\begin{equation}\label{eq:nint}
\nint T:=\mathrm{Res}_{s=0}\tr(T|D|^{-s})
\end{equation}
is tracial on $\Psi^0$. We also recall the definition of `smoothing operators' $\op$, 
$$
\op:=\{T\in\opz\,|\,|D|^kT\in\opz\,,\;\forall\;k\in\N\}\;.
$$
The class $\op$ is a two-sided $*$-ideal in the $*$-algebra $\opz$, is $\delta$-invariant
and then in the smooth domain of $\delta$. If $T$ is a smoothing operator, $\zeta_T(s)$
is holomorphic on $\C$ and (\ref{eq:nint}) vanishes.
Thus, elements in $\op$ can be neglected when computing the dimension spectrum
and residue.
Finally, we note that if the metric dimension is finite, rapid decay matrices -- in a basis
of eigenvectors for $D$ with eigenvalues in increasing order -- are smoothing operators.

\bigskip

In analogy with the notion of spin manifold, one asks for the existence of a real structure
$J$ on a spectral triple $(\A,\HH,D)$. Motivated by the examples of real spectral triples
on Podle{\'s} spheres~\cite{DLPS,DDLW} and on $SU_q(2)$ \cite{DLSvSV}, we use the following
weakened definition of real structure.

\begin{df}
A \emph{real structure} is an antilinear isometry $J$ on $\HH$ such that
$\forall\;a,b\in\A$,
$$
J^2=\pm 1\;,\qquad
JD=\pm DJ\;,\qquad
[a,JbJ^{-1}]\subset\mathcal{I}\;,\qquad
[[D,a],JbJ^{-1}]\subset\mathcal{I}\;.
$$
If the spectral triple is even with grading $\gamma$, we impose the further
relation $J\gamma=\pm\gamma J$.
\end{df}

\noindent
The signs `$\pm$' are determined by the dimension of the geometry~\cite{ConnesReal}.
A real spectral triple of dimension $4$ corresponds to the choices $J^2=-1$, $JD=DJ$
and $J\gamma=\gamma J$.

\noindent
The set $\mathcal{I}$ is a suitable two-sided ideal in the algebra $\opz$ of
`order zero' operators which is made of `infinitesimals'.
The original definition~\cite{ConnesReal} corresponds to $\mathcal{I}=0$;
while in examples coming from quantum groups~\cite{DLPS,DLSvSV,DDLW}
one usually takes $\mathcal{I}=\op$.

Let $F:=D|D|^{-1}$ be the sign of $D$; if $(\A,\HH,D)$ is a regular even spectral triple, the datum
$(\A,\HH,F,\gamma)$ is an even Fredholm module.
We say that the Fredholm module is \emph{$p$-summable} if $p\geq 1$ and,
for all $a\in\A$, $[F,a]$ belongs to the $p$-th Schatten-von Neumann ideal
$\mathcal{L}^p(\HH)$ of compact operators $T$ such that $|T|^p$ is of trace class. Associated with a $p$-summable even Fredholm module there are cyclic cocycles defined by
\begin{equation}\label{eq:ch}
\mathrm{ch}^F_n(a_0,\ldots,a_n)=\tfrac{1}{2\,n!}\,\Gamma(\tfrac{n}{2}+1)\,\tr(\gamma F[F,a_0]\ldots[F,a_n])\;,
\end{equation}
for all even integers $n\geq p-1$. By composing it with a matrix trace, 
$\mathrm{ch}^F_n$ is canonically extended to matrices with entries in $\A$.
The pairing with elements $[e]\in K_0(\A)$, given by $\mathrm{ch}^F_n(e,e,\ldots,e)$ build up to an integer-valued map $\mathrm{ch}^F\!([e])$ which depends only on
the class $[e]$ and which yields the index of the
Dirac operator $D$ twisted with the projection $e$ (for further details see~\cite{CBook}).

\bigskip

Finally, we turn now to symmetries; these 
will be implemented by an action of a Hopf $*$-algebra.  
Firstly, let $\mathcal{V}$ be a dense linear subspace of a Hilbert space $\HH$ with inner product 
$\inner{\,,\,}$, and let $\U$ be a $*$-algebra.
An (unbounded) $*$-representation of $\U$ on $\mathcal{V}$ is a homomorphism 
$\lambda:\U\to\mathrm{End}(\mathcal{V})$ such that $\inner{\lambda(h)v,w}=\inner{v,\lambda(h^*)w}$ for all
$v,w\in\mathcal{V}$ and all $h\in\U$. From now on, the symbol $\lambda$ will be omitted. Next, let $\U=(\U,\Delta,\varepsilon,S)$ be a Hopf $*$-algebra and let $\A$ be a left $\U$-module
$*$-algebra, i.e., there is a left action $\az$ of $\U$ on $\A$ satisfying 
$$
h\az ab=(h_{(1)}\az a)(h_{(2)}\az b)\;,\qquad
h\az 1=\varepsilon(h)1\;,\qquad
h\az a^*=\{S(h)^*\az a\}^*\;,
$$
for all $h\in\U$ and $a,b\in\A$. As customary, $\Delta(h)=h_{(1)}\otimes h_{(2)}$.

A $*$-representation of $\A$ on $\mathcal{V}$ is called $\U$-equivariant if there exists a
$*$-representation of $\U$ on $\mathcal{V}$ such that, for all $h\in\U$, $a\in\A$ and
$v\in\mathcal{V}$, it happens that
$$
hav=(h_{(1)}\az a)\,h_{(2)}v \;.
$$
Given $\U$ and $\A$ as above, the left crossed product $*$-algebra $\A\rtimes\U$ is
defined as the $*$-algebra generated by the two $*$-subalgebras $\A$ and $\U$ with
crossed commutation relations
$$
ha=(h_{(1)}\az a)h_{(2)}\;,\quad\forall\;h\in\U,\;a\in\A\,.
$$
Thus, $\U$-equivariant $*$-representations of $\A$ correspond to
$*$-representations of $\A\rtimes\U$.

A linear operator $D$ defined on $\mathcal{V}$ is said to be equivariant if 
it commutes with $\U$, i.e.,  
\begin{equation}
Dhv=hDv
\end{equation} 
for all $h\in\U$ and $v\in\mathcal{V}$. On the other hand, an antilinear operator $T$
defined on $\mathcal{V}$ is called equivariant if it satisfies the relation
\begin{equation}
Thv=S(h)^*Tv \,\, , 
\end{equation} 
for all $h\in\U$ and $v\in\mathcal{V}$, where $S$ denotes the antipode of $\U$.
Notice that if $T$ is an equivariant antilinear operator, its square $T^2$ is an equivariant
linear operator, but $T^*T$ is not an equivariant linear operator unless $S^2=1$.

We use all these equivariance requirements in the following definition (see also \cite{Sitarz}).

\begin{df}
Let $\U$ be a Hopf $*$-algebra and $\A$ a left $\U$-module $*$-algebra.
A (real, even) spectral triple $(\A,\HH,D,\gamma,J)$ is called \emph{equivariant} if 
$\U$ is represented on a dense subspace $\mathcal{V}$ of $\HH$, $\mathcal{V}\subset\mathrm{dom}\,D$,
the representation of $\U$ commutes with the grading $\gamma$,
the restriction of the representation of $\A$ on $\mathcal{V}$ is $\U$-equivariant,
the operator $D$ is equivariant  and $J$ is the antiunitary part
of the polar decomposition of an equivariant antilinear operator.
\end{df}

\subsection{Projective module description of equivariant re\-pre\-sen\-ta\-tions}\label{app:B}

In order to construct the analogues of the modules of chiral spinors on the sphere
$S^4_q$ we need some properties of finitely generated projective modules over algebras
having quantum group symmetries.

Let $\mU$ be a Hopf $*$-algebra, $\A$ be an $\mU$-module $*$-algebra and $\varphi:\A\to\C$ be
an invariant faithful state (i.e.~$\varphi$ is linear, $\varphi(a^*a)>0$ for all
nonzero $a\in\A$, and $\varphi(h\az a)=\epsilon(h)\varphi(a)\;\forall\;a\in\A$
and $h\in\mU$). Suppose also that there exists $\kappa\in\mathrm{Aut}(\A)$
such that the `twisted' cyclicity
$$
\varphi(ab)=\varphi\bigl(b\,\kappa(a)\bigr)
$$
holds for all $a,b\in\A$. Instances of this situation are provided by 
subalgebras of compact quantum group algebras with $\varphi$ the Haar state and $\kappa$ the modular involution. KMS states in Thermal Quantum Field Theory provide additional examples.
In particular, for the case $\A=\A(S^4_q)$ and $\mU=\Uq$, $\varphi$
comes from the Haar functional of $\A(SO_{q^2}(5))$ and the modular
automorphism is $\kappa(a)=K_1^8K_2^6\az a$~\cite[Sect.~11.3.4]{Kli}.

For $N\in\N$, let $\A^N:=\A\otimes\C^N$ be the linear space with
elements $v=(v_1,\ldots,v_N)$, $v_i\in\A$, and $\C$-valued inner product given by
\begin{equation}\label{eq:in}
\inner{v,w}:=\sum\nolimits_{i=1}^N\varphi(v_i^*w_i)\;.
\end{equation}

\begin{lemma}\label{lemma:ip}
Let $\sigma:\U\to\mathrm{Mat}_N(\C)$ be a $*$-representation.
The formul{\ae}:
\begin{equation}\label{eq:starrep}
(a.v)_i:=av_i\;,\qquad
(h.v)_i:=\sum\nolimits_{j=1}^N(h_{(1)}\az v_j)\sigma_{ij}(h_{(2)})\;,
\end{equation}
for all $a,v\in\A$ and $h\in\U$ (and $i=1,\dots,N$), define a $*$-representation of the crossed product algebra
$\A\rtimes\U$ on the linear space $\A^N$.
\end{lemma}
\begin{prova}
The inner product allows us to define the adjoint of an element of $\A\rtimes\U$
in the representation on $\A^N$. For $x\in\mathrm{End}(\A^N)$,  its adjoint denoted with $x^\dag$, is defined 
$$
\inner{\smash[t]{x^\dag.v},w}:=\inner{v,x.w}\;, \qquad \forall \;, v,w\in\A^N\;.
$$
Recall that being a $*$-representation means that $x^\dag.v=x^*.v$ for any
operator $x$ and any $v\in\A^N$.

The nontrivial part of the proof consists in showing that $h^\dag.v=h^*.v$ for all
$h\in\U$ and $v\in\A$. For $N>1$ we are considering the Hopf tensor product of the
$N=1$ representation with a matrix representation that is a $*$-representation by
hypothesis. Thus, it is enough to take $N=1$.

The $\U$-invariance of $\varphi$ implies:
$$
\epsilon(h)\inner{v,w}=\varphi\bigl(h\az(v^*w)\bigr)=
\varphi\bigl((h_{(1)}\az v^*)(h_{(2)}\az w)\bigr)\;.
$$
But $h_{(1)}\az v^*=\{S(h_{(1)})^*\az v\}^*$ by definition of module
$*$-algebra. Then,
$$
\epsilon(h)\inner{v,w}=\inner{S(h_{(1)})^*.v,h_{(2)}.w}=
\inner{v,S(h_{(1)})^{*\dag}h_{(2)}.w}\;.
$$
We deduce that for all $h\in\U$ one has that
\begin{equation}\label{eq:newS}
S(h_{(1)})^{*\dag}h_{(2)}=\epsilon(h)\;.
\end{equation}
Recall that the convolution product `$\star$' for any $F,G\in\mathrm{End}(\U)$ 
is defined by
$$
(F\star G)(h):=F(h_{(1)})G(h_{(2)})\qquad\forall\;h\in\U \;;
$$
and $(\mathrm{End}(\U),\star)$ is an associative algebra with unity given
by the endomorphism $h\mapsto\epsilon(h)1_{\U}$,  with $S$  a
left and right inverse for $id_{\U}$ in $(\mathrm{End}(\U),\star)$, that is
$$
S\star id_{\U}=1_{\U}\epsilon=id_{\U}\star S\;.
$$
Let $S'\in\mathrm{End}(\U)$ be the composition $S':=\dag\circ *\circ S$.
Equation \eqref{eq:newS} implies that $S'$ is a left inverse for $id_{\U}$:
$$
S'\star id_{\U}=1_{\U}\epsilon\;.
$$
Applying $\star S$ to the right of both members of this equation and using
$id_{\U}\star S=1_{\U}\epsilon$ we get $S'=S$ as endomorphisms of $\U$,
i.e.~$S(h)^{*\dag}=S(h)$ for all $h\in\U$.

Now, the antipode of a Hopf $*$-algebra is invertible,
with $S^{-1}=*\circ S\circ *$, thus we arrive at $h^{*\dag}=h$ for all $h\in\U$.
Replacing $h$ with $h^*$ we prove that $h^\dag=h^*$ for all $h\in\U$, and
this concludes the proof.
\end{prova}

Now, let $e=(e_{ij})\in\mathrm{Mat}_N(\A)$ be an $N\times N$ matrix
with entries $e_{ij}\in\A$. Let $\pi:\A^N\to\A^N$ be the (linear)
endomorphism defined by:
\begin{equation}\label{eq:pi}
\pi(v)_j:=\sum\nolimits_{i=1}^Nv_ie_{ij}\;,
\end{equation}
for all $v\in\A^N$ and $j=1,\dots,N$.
Since $\A$ is associative, left and right multiplication commute
and $\pi(av)=a\pi(v)$ for all $a\in\A$ and $v\in\A^N$.
Thus we have the following lemma.

\begin{lemma}
The map $\pi$ defined by \eqref{eq:pi} is an $\A$-module map.
\end{lemma}

Recall that an endomorphism $p$ of an inner
product space $V$ is a projection (not necessarily orthogonal)
if $p\circ p=p$. A projection $p$ is \emph{orthogonal}
if the image of $p$ and $id_V-p$ are orthogonal with respect
to the inner product of $V$, and this happens exactly when $p^\dag=p$.

A simple computation shows that the map $\pi$ in 
\eqref{eq:pi} is a projection if{}f $e^2=e$, that is the matrix $e\in\mathrm{Mat}_N(\A)$ is an idempotent. Now we use
the twisted-cyclicity of $\varphi$ to deduce:
$$
\inner{v,\pi^\dag(w)}=\inner{\pi(v),w}=
\sum\nolimits_{ij}\varphi(e_{ij}^*v_i^*w_j)=
\sum\nolimits_{ij}\varphi\bigl(v_i^*w_j\kappa(e_{ij}^*)\bigr) \;, 
$$
for all $v,w\in\A^N$. Hence the adjoint $\pi^\dag$ of the endomorphism $\pi$
is given by
$$
\pi^\dag(w)_i=\sum\nolimits_{j=1}^Nw_j\kappa(e_{ij}^*)\;.
$$
Let $e^*$ be the matrix with entries $(e^*)_{jk}:=e_{kj}^*$.
We have the following lemma.

\begin{lemma}\label{lemma:next}
The endomorphism $\pi$ in \eqref{eq:pi} is an orthogonal
projection if{}f $\,e^2=e=\kappa(e^*)$.
\end{lemma}

Next, we  determine a sufficient condition for the endomorphism $\pi$ to be
not only an $\A$-module map, but also an $\U$-module map.

\begin{lemma}
With `$\phantom{|}^t$' denoting transposition, if
\begin{equation}\label{eq:cov}
h\az e=\sigma(h_{(1)})^t\,e\,\sigma(S^{-1}(h_{(2)}))^t \;, 
\end{equation}
for all $h\in\U$, the endomorphism $\pi$ in \eqref{eq:pi} is an $\U$-module map.
\end{lemma}
\begin{prova}
Equation (\ref{eq:cov}) can be rewritten as,
$$
h\az e_{ij}=\sum\nolimits_{kl}\sigma_{ki}(h_{(1)})\,e_{kl}
\,\sigma_{jl}(S^{-1}(h_{(2)}))\;;
$$
by using it into the definition \eqref{eq:starrep}
one checks that $\pi(h.v)=h.\pi(v)$ for all $h\in\U$ and $v\in\A^N$.
\end{prova}

When Lemma \ref{lemma:next} and Eq.~(\ref{eq:cov}) are satisfied,
the orthogonal projections $\pi$ and \linebreak $\pi^\perp=1-\pi\,$
split $\A^N$ into the orthogonal sum of two sub $*$-representations
$\pi(\A^N)$ and $\pi^\perp(\A^N)$ of $\A\rtimes\U$.
The next lemma gives a (quite obvious) sufficient condition for $\pi(\A^N)$
and $\pi^\perp(\A^N)$ to be not equivalent as representations of $\A$.
Recall that an isomorphism of $\A$-modules is an invertible $\A$-linear
map, so isomorphic modules correspond to equivalent representations.

\begin{lemma}\label{lemma:inv}
Let $(\A,\HH,F,\gamma)$ be an even Fredholm module over $\A$.
If $\mathrm{ch}^F\!([e])\neq 0$,
the $\A$-modules $\pi(\A^N)$ and $\pi^\perp(\A^N)$
are not equivalent.
\end{lemma}
\begin{prova}
The map $K_0(\A)\to\Z$, $[e]\mapsto\mathrm{ch}^F\!([e])$ is an
homomorphism.
Suppose $\pi(\A^N)$ and $\pi^\perp(\A^N)$ are isomorphic
$\A$-modules, then $[e]=[1-e]$ and
$\mathrm{ch}^F\!([1-e])=\mathrm{ch}^F\!([e])$.

But from Eq.~(\ref{eq:ch}),
$\mathrm{ch}^F\!([1-e])=-\mathrm{ch}^F\!([e])$ (since $[F,1-e]=-[F,e]$
and $n$ is even).

Hence $\mathrm{ch}^F\!([e])=0$, and this concludes the proof by contradiction.
\end{prova}


\section{The symmetry Hopf algebra $\Uq$}\label{sec:rep}
Let $0<q<1$. We call $\Uq$ the real form of the Drinfeld-Jimbo deformation
of $so(5,\C)$, corresponding to the Euclidean signature $(+,+,+,+,+)$;
it is a real form of the Hopf algebra called $\breve{U}_q(so(5,\C))$
in~\cite[Sect.~6.1.2]{Kli}. As a $*$-algebra, 
$\Uq$ is generated by $\{K_i=K_i^*,K_i^{-1},E_i,F_i:=E_i^*\}_{i=1,2}$ ($i\to 3-i$ with
respect to the notations of~\cite{Kli}), with relations:
\begin{equation*}\begin{array}{c}
[K_1,K_2]=0\;,\qquad
K_iK_i^{-1}=K_i^{-1}K_i=1\;,
\\ \rule{0pt}{20pt}
[E_i,F_j]=\delta_{ij}\frac{K_j^2-K_j^{-2}}{q^j-q^{-j}}\;,
\\ \rule{0pt}{20pt}
K_iE_iK_i^{-1}=q^iE_i\;,\qquad
K_iE_jK_i^{-1}=q^{-1}E_j\;\;\mathrm{if}\;i\neq j\;,
\end{array}\end{equation*}
together with the ones obtained by conjugation and Serre relations, 
explicitly, given by
\begin{subequations}\label{eq:Serre}
\begin{align}
E_1E_2^2-(q^2+q^{-2})E_2E_1E_2+E_2^2E_1 &=0\;, \\*
E_1^3E_2-(q^2+1+q^{-2})(E_1^2E_2E_1-E_1E_2E_1^2)-E_2E_1^3 &=0\;,
\end{align}
\end{subequations}
together with their adjoints. These relations can be written in a more compact
form by defining $[a,b]_q:=q^2ab-ba$. Then, (\ref{eq:Serre}) are equivalent to
$$
[E_2,[E_1,E_2]_q]_q=0\;,\qquad
[E_1,[E_1,[E_2,E_1]_q]_q]=0\;.
$$
The Hopf algebra structure $(\Delta,\epsilon,S)$ of $\Uq$ is given by:
\begin{equation*}
\begin{array}{c}
\Delta K_i=K_i\otimes K_i\;,\quad \Delta E_i=E_i\otimes K_i+K_i^{-1}\otimes E_i\;, 
\\ \rule{0pt}{3.5ex}\epsilon(K_i)=1\;,\quad \epsilon(E_i)=0\;,\\ \rule{0pt}{3.5ex}
S(K_i)=K_i^{-1}\;,\quad S(E_i)=-q^iE_i\;.
\end{array}
\end{equation*}
For each non negative $n_1,n_2$ such that $n_2\in\frac{1}{2}\N$ and $n_2-n_1\in\N$ there
is an irreducible representation of $\Uq$ whose representation space we denote
$V_{(n_1,n_2)}$.
We call it ``the representation with highest weight $(n_1,n_2)$'' since
the highest weight vector is an eigenvector of $K_1$ and $K_1K_2$
with eigenvalues $q^{n_1}$ and $q^{n_2}$, respectively.

Irreducible representations with highest weight $(0,l)$ and $(\frac{1}{2},l)$
(the ones that we need explicitly) can be found in~\cite{Cha} and are recalled 
presently. 
Let us use the shorthand notation $V_l:=V_{(0,l)}$ if $l\in\N$ and $V_l:=V_{(\frac{1}{2},l)}$
if $l\in\N+\frac{1}{2}$.
The vector space $V_l$, for all $l\in\frac{1}{2}\N$, has orthonormal basis $\ket{l,m_1,m_2;j}$,
where the labels $(j,m_1,m_2)$ satisfy the following constraints.
For $l\in\N$:
$$
j=0,1,\ldots,l\;,\qquad
j-|m_1|\in\N\;,\qquad
l-j-|m_2|\in 2\N\;,
$$
while for $l\in\N+\frac{1}{2}$:
$$
j=\tfrac{1}{2},\tfrac{3}{2},\ldots,l-1,l\;,\qquad
j-|m_1|\in\N\;,\qquad
l+\tfrac{1}{2}-j-|m_2|\in\N\;.
$$
Notice that for any admissible $(l,m_1,m_2,j)$ there exists a unique $\epsilon\in\{0,\pm\frac{1}{2}\}$
such that \mbox{$l+\epsilon-j-m_2\in 2\N$} (that is, $\epsilon=0$ if $l\in\N$ and
$\epsilon=\frac{1}{2}(-1)^{l+\frac{1}{2}-j-m_2}$ if $l\in\N+\frac{1}{2}$). We shall need the
coefficients,
\begin{subequations}
\begin{align}
a_l(j,m_2) &=\frac{1}{[2]}\sqrt{\frac{[l-j-m_2+\epsilon][l+j+m_2+3+\epsilon]}
  {[2(j+|\epsilon|)+1][2(j-|\epsilon|)+3]}} \;, \\
b_l(j,m_2) &=2|\epsilon|\,\frac{\sqrt{[l-\epsilon(2j+1)-m_2+1][l-\epsilon(2j+1)+m_2+2]}}{[2j][2j+2]}
 \;,\label{eq:abc} \\
c_l(j,m_2) &=\frac{(-1)^{2\epsilon}}{[2]}\sqrt{\frac{[l-j+m_2+2-\epsilon][l+j-m_2+1-\epsilon]}
  {[2(j+|\epsilon|)-1][2(j-|\epsilon|)+1]}} \;,
\end{align}
\end{subequations}
where, as usual, $[z]:=(q^z-q^{-z})/(q-q^{-1})$ denotes the $q$-analogue of $z\in\C$.

The $*$-representation $\sigma_l:\Uq\to\mathrm{End}(V_l)$ is defined by the rules,
\begin{align*}
\sigma_l(K_1)\ket{l,m_1,m_2;j} &=q^{m_1}\ket{l,m_1,m_2;j} \;, \\
\sigma_l(K_2)\ket{l,m_1,m_2;j} &=q^{m_2-m_1}\ket{l,m_1,m_2;j} \;, \\
\sigma_l(E_1)\ket{l,m_1,m_2;j} &=\sqrt{[j-m_1][j+m_1+1]}\ket{l,m_1+1,m_2;j} \;, \\
\sigma_l(E_2)\ket{l,m_1,m_2;j} &=
\sqrt{[j-m_1+1][j-m_1+2]}\,a_l(j,m_2)\ket{l,m_1-1,m_2+1;j+1} \\ &+
\sqrt{[j+m_1][j-m_1+1]}\,b_l(j,m_2)\ket{l,m_1-1,m_2+1;j} \\ &+
\sqrt{[j+m_1][j+m_1-1]}\,c_l(j,m_2)\ket{l,m_1-1,m_2+1;j-1} \;.
\end{align*}
When there is no risk of ambiguity the representation symbol $\sigma_l$ will be suppressed.

For $l\in\N$ the representation $\sigma_l$ is \emph{real}. That is, there is an
antilinear map $C:V_l\to V_l$, which satisfies $C^2=1$ and $C\sigma_l(h)C=\sigma_l(S(h)^*)$.
This map is explicitly given by 
\begin{equation}\label{eq:C}
C\ket{l,m_1,m_2;j}:=(-q)^{m_1}q^{3m_2}\ket{l,-m_1,-m_2;j}  \;.
\end{equation}

The operator
\begin{equation}\label{eq:Cdef}
\mathcal{C}_1:=q^{-1}K_1^2+qK_1^{-2}+(q-q^{-1})^2E_1F_1\;,
\end{equation}
is a Casimir for the subalgebra generated by $(K_1,K_1^{-1},E_1,F_1)$.
For future reference, we note the action of $\mathcal{C}_1$ on
a vector of $V_l$, with $l\in\frac{1}{2}\N$; it is 
\begin{equation}\label{eq:Caz}
\mathcal{C}_1\ket{l,m_1,m_2;j}=(q^{2j+1}+q^{-2j-1})\ket{l,m_1,m_2;j} \;.
\end{equation}

\section{The orthogonal quantum  $4$-Sphere}\label{sec:S4q}
\begin{df}[\cite{FRT}]\label{df:S4q}
We call \emph{orthogonal quantum  $4$-sphere} the virtual space underlying the algebra
$\A(S^4_q)$ generated by $x_0=x_0^*,x_i$ and $x_i^*$ (with $i=1,2$), with commutation
relations:
\begin{align*}
&x_ix_j     = q^2x_jx_i   \;, \qquad\quad\forall\;\;0\leq i<j\leq 2 \;, \\
&x_i^*x_j   =q^2x_jx_i^*  \;, \qquad\quad\forall\;\;i\neq j \;, \\
&[x_1^*,x_1]=(1-q^4)x_0^2 \;, \\
&[x_2^*,x_2]=x_1^*x_1-q^4x_1x_1^* \;, \\
&x_0^2+x_1x_1^*+x_2x_2^*=1 \;.
\end{align*}
\end{df}

\noindent
The original notations of Fadeev-Reshetikhin-Takhtadzhyan~\cite[Eq.~(1.14)]{FRT}
can be obtained by defining $x'_1:=x_2^*$, $x'_2:=x_1^*$, $x'_3:=\sqrt{q(1+q^2)}\,x_0$,
$x'_4:=x_1$, $x'_5:=x_2$ and $q':=q^2$.
The notations in~\cite[Eq.~(2.1)]{LH} can be obtained by the replacement
$x_i\mapsto x_i^*$ and $q^2\mapsto q^{-1}$.

\bigskip

In the next propositions we summarize some well known facts.

\begin{prop}
The algebra $\A(S^4_q)$ is an $\Uq$-module $*$-algebra for the action
given by:
\begin{equation*}\begin{array}{c}
K_i\az x_i=qx_i \;, i=1,2\;, \\
\rule{0pt}{18pt}
K_2\az x_1=q^{-1}x_1\;,  \\
\rule{0pt}{18pt}
E_1\az x_0=q^{-1/2}x_1\;,\qquad
E_2\az x_1=x_2\;,  \\
\rule{0pt}{18pt}
F_1\az x_1=q^{1/2}[2]x_0 \;,\quad
F_1\az x_0=-q^{-3/2}x_1^*\;\;\quad
F_2\az x_2=x_1\;,
\end{array}\end{equation*}
while $K_i\az x_j=x_j$, $E_i\az x_j=0$ and $F_i\az x_j=0$ in all other cases.
\end{prop}

\noindent
Notice that the action on the $x_i^*$'s is determined by compatibility with the involution:
$$
K_i\az a^*=\{K_i^{-1}\az a\}^* \;,\qquad
E_1\az a^*=\{-qF_1\az a\}^* \;,\qquad
E_2\az a^*=\{-q^2F_2\az a\}^* \;.
$$

\begin{prova}
The bijective linear map from the linear span of $\{x_i,x_i^*\}$
to the representation space $V_1$ defined (modulo a global
proportionality constant) by 
$$
x_2\mapsto\ket{0,1;0}\,,\;
x_1\mapsto\ket{1,0;1}\,,\;
x_0\mapsto (q[2])^{-1/2}\ket{0,0;1}\,,\;
x_1^*\mapsto -q\ket{-1,0;1}\,,\;
x_2^*\mapsto q^3\ket{0,-1;0}\,,
$$
is a unitary equivalence of $\Uq$-modules (here unitary means that
the real structure $C$ on $V_1$ is implemented by the $*$ operation
on $x_i$'s). This guarantees that the free $*$-algebra $\C\!\inner{x_i,x_i^*}$
generated by $\{x_i,x_i^*\}$ is an $\Uq$-module $*$-algebra.

The degree $\leq 2$ polynomials generating the ideal which
defines $\A(S^4_q)$ span the real representations $V_0$ and
$V_{(1,1)}$, inside the tensor product $V_1\otimes V_1$.
The quotient $*$-algebra of $\C\!\inner{x_i,x_i^*}$ by this ideal,
$\A(S^4_q)$, is then an $\Uq$-module $*$-algebra.
\end{prova}

\begin{prop}\label{prop:S4q}
There is an isomorphism $\,\A(S^4_q)\simeq\bigoplus_{l\in\N}V_l\,$ of $\Uq$ left modules.
\end{prop}
\begin{prova}
A linear basis for $\A(S^4_q)$ is made of monomials $x_0^{n_0}x_1^{n_1}(x_1^*)^{n_2}x_2^{n_3}$
with $n_0,n_1,n_2\in\N$, $n_3\in\Z$ and with the notation $x_2^{n_3}:=(x_2^*)^{|n_3|}$
if $n_3<0$. Using this basis one proves that a weight vector of $\A(S^4_q)$
is annihilated by both $E_1$ and $E_2$ if and only if it is of the
form $x_2^l$, $l\in\N$. Thus, highest weight vectors are proportional
to $x_2^l$ and the algebra decomposes as multiplicity free direct sum
of highest weight representations with weights $(0,l)$.
\end{prova}

The algebra $\A(S^4_q)$ has two inequivalent irreducible infinite dimensional representations. The representation space is the Hilbert space $\ell^2(\N^2)$ and the representations are given by  
\begin{align}\label{4irreps}
x_0\ket{k_1,k_2}_\pm &:=\pm q^{2(k_1+k_2)}\ket{k_1,k_2}_\pm \;, \nonumber \\
x_1\ket{k_1,k_2}_\pm &:=q^{2k_2}\sqrt{1-q^{4(k_1+1)}}\ket{k_1+1,k_2}_\pm \;,\\
x_2\ket{k_1,k_2}_\pm &:=\sqrt{1-q^{4(k_2+1)}}\ket{k_1,k_2+1}_\pm \nonumber \;.
\end{align}
The direct sum of these representations, with obvious grading $\gamma$ and  
operator $F$ given by $F\ket{k_1,k_2}_\pm:=\ket{k_1,k_2}_\mp$, constitutes
a $1$-summable Fredholm module over $\A(S^4_q)$. 

\bigskip

In the sequel we shall need both the quantum space $SU_q(2)$ as well as
the equatorial Podle\'s sphere, whose algebras are given in~\cite{Wor}
and~\cite{Pod} respectively.

\begin{df}
The algebra $\A(SU_q(2))$ of polynomial  functions on $SU_q(2)$ is
the $*$-algebra generated by $\alpha,\beta$ and their
adjoints, with relations:
$$
\beta\alpha=q\alpha\beta \;,\quad
\beta^*\alpha=q\alpha\beta^* \;,\quad
[\beta,\beta^*]=0 \;,\quad
\alpha\alpha^*+\beta\beta^*=1 \;,\quad
\alpha^*\alpha+q^2\beta^*\beta=1 \;.
$$
We call \emph{equatorial Podle\'s sphere} the virtual space underlying
the $*$-algebra $\A(S^2_q)$ generated by $A=A^*$, $B$ and $B^*$ with
relations:
$$
AB=q^2BA \;,\qquad
BB^*+A^2=1 \;,\qquad
B^*B+q^4A^2=1 \;.
$$
\end{df}

\begin{prop}
There is a $*$-algebra morphism $\varphi:\A(S^4_q)\to\A(SU_q(2))\otimes\A(S^2_q)$
defined by:
\begin{align}
\varphi(x_0) &=-(\alpha\beta+\beta^*\alpha^*)\otimes A\;, \nonumber \\
\varphi(x_1) &=\bigl(-\alpha^2+q\,(\beta^*)^2\bigr)\otimes A\;, \label{eq:embed} \\
\varphi(x_2) &=1\otimes B\;. \nonumber
\end{align}
\end{prop}
\begin{prova}
One proves by direct computation that the five elements $\varphi(x_i),\varphi(x_i)^*$ satisfy
all the defining relations of $\A(S^4_q)$.
\end{prova}


\section{The modules of chiral spinors}\label{sec:cs}

We apply the general theory of Sect.~\ref{app:B}, 
to the case $\A=\A(S^4_q)$ and $\U=\Uq$. Recall that in this case
$\kappa(a)=K_1^8K_2^6\az a$ is the modular automorphism.
We shall use the notations of Sect.~\ref{sec:rep} for the irreducible
representations $(V_l,\sigma_l)$ of $\Uq$. 

By Proposition \ref{prop:S4q} we have the equivalence
\mbox{$\A(S^4_q)\simeq\bigoplus_{l\in\N}V_l$} as left $\Uq$-modules.
Using Lemma~\ref{lemma:ip} for $N=1$, we deduce that on the vector space
$\bigoplus_{l\in\N}V_l$ there exists at least one $*$-representation of
the crossed product $\A(S^4_q)\rtimes\Uq$ that extends the $*$-representation
$\,\bigoplus_{l\in\N}\sigma_l\,$ of $\Uq$.

Let $e\in\mathrm{Mat}_4(\A(S^4_q))$ be the following idempotent:
\begin{equation}\label{eq:P}
e:=\frac{1}{2}\left(\begin{array}{cccc}
   1+x_0       &  q^3x_2   & -qx_1        & 0 \\
   q^{-3}x_2^* &  1-q^2x_0 &  0           & q^3x_1 \\
  -q^{-1}x_1^* &  0        &  1-q^2x_0    & q^3x_2 \\
   0           &  qx_1^*   &  q^{-3}x_2^* & 1+q^4x_0
  \end{array}\right)\;.
\end{equation}
By direct computation one proves that $K_1^8K_2^6\az e^*=e=e^2$ and then,
by Lemma \ref{lemma:next},
$e$ defines an orthogonal projection $\pi$, by equation~\eqref{eq:pi}, on the linear space
$\A(S^4_q)^4$ with inner product~\eqref{eq:in}.
Next, let $\sigma:\Uq\to\mathrm{Mat}_4(\C)$ be the $*$-representation defined by
\begin{subequations}\label{eq:spin}
\begin{equation}
\sigma(K_1)=\ma{
  q^{1/2} & 0 & 0 & 0 \\
  0 & q^{1/2} & 0 & 0 \\
  0 & 0 & q^{-1/2} & 0 \\
  0 & 0 & 0 & q^{-1/2}
}\;,\qquad
\sigma(K_2)=\ma{
  1 & 0 & 0 & 0 \\
  0 & q^{-1} & 0 & 0 \\
  0 & 0 & q & 0 \\
  0 & 0 & 0 & 1
}\;,
\end{equation}
\begin{equation}
\sigma(E_1)=\ma{
  0 & 0 & 1 & 0 \\
  0 & 0 & 0 & 1 \\
  0 & 0 & 0 & 0 \\
  0 & 0 & 0 & 0
}\;,\qquad
\sigma(E_2)=\ma{
  0 & 0 & 0 & 0 \\
  0 & 0 & 0 & 0 \\
  0 & 1 & 0 & 0 \\
  0 & 0 & 0 & 0
}\;.
\end{equation}
\end{subequations}
Again, by direct computation one proves that:
\begin{subequations}\label{eq:sub}
\begin{align}
K_i\az e &=\sigma(K_i)\,e\,\sigma(K_i)^{-1} \;,\\
E_i\az e &=\sigma(F_i)\,e\,\sigma(K_i)^{-1}-q^{-i}\sigma(K_i)^{-1}e\,\sigma(F_i) \;.
\end{align}
\end{subequations}
Since $\sigma(K_i)=\sigma(K_i)^t$ and $\sigma(F_i)=\sigma(E_i)^t$, we conclude
that condition (\ref{eq:cov}) is satisfied and that $\pi$ and $\pi^\perp=1-\pi$ project
$\A(S^4_q)^4$ onto sub $*$-representations of $\A(S^4_q)\rtimes\Uq$.

We state the main proposition of this section.

\begin{prop}\label{prop:main}
There exists two inequivalent representations of the crossed product
algebra $\A(S^4_q)\rtimes\Uq$ on $\bigoplus_{l\in\N+\frac{1}{2}}V_l$
that extend the representation $\,\bigoplus_{l\in\N+\frac{1}{2}}\sigma_l\,$
of $\Uq$.
\end{prop}

The proof is in two steps. We first prove (in Lemma \ref{lemma:ptwo})
that $\pi(\A(S^4_q)^4)$ and $\pi^\perp(\A(S^4_q)^4)$ are not equivalent
as representations of the algebra $\A(S^4_q)$. Then we prove (in Lemma
\ref{lemma:pone}) that as $\Uq$ representations they are both equivalent
to $\bigoplus_{l\in\N+\frac{1}{2}}V_l$.

\begin{lemma}\label{lemma:ptwo}
The idempotent $e$ in (\ref{eq:P}) splits $\A(S^4_q)^4$ into two \emph{inequivalent}
$*$-representations of the crossed product algebra $\A(S^4_q)\rtimes\Uq$.
\end{lemma}

\begin{prova}
To prove the statement we apply Lemma~\ref{lemma:inv}.
We use the Fredholm module associated to the representation on
$\ell^2(\N)\oplus\ell^2(\N)$ defined by Eq.~(\ref{4irreps}).
One has
\begin{align*}
\mathrm{ch}^F\!([e])
&=\tfrac{1}{2}\tr_{\ell^2(\N)\otimes\C^8}(\gamma F[F,e]) \\
\rule{0pt}{14pt}
&=\tfrac{1}{4}(1-q^2)^2\tr_{\ell^2(\N)\otimes\C^2}(\gamma F[F,x_0]) \\
\rule{0pt}{14pt}
&=(1-q^2)^2\sum\nolimits_{k_1,k_2\in\N}q^{2(k_1+k_2)}=1\;.
\end{align*}
The statement of Proposition~\ref{prop:main} follows from the obvious observation
that if the two representations of the crossed product algebra were
equivalent, their restrictions to representations of $\A(S^4_q)$
would be equivalent too.
\end{prova}

\begin{lemma}\label{lemma:pone}
$\pi(\A(S^4_q)^4)\simeq\pi^\perp(\A(S^4_q)^4)\simeq\bigoplus_{l\in\N+\frac{1}{2}}V_l$
as $\Uq$ representations.
\end{lemma}
\begin{prova}
In this proof, `$\simeq$' always means equivalence of representations of $\Uq$.

Since $\sigma$ in \eqref{eq:spin} is unitary equivalent to the spin representation $V_{1/2}$,
the representation of $\Uq$ on $\A(S^4_q)^4$ is the Hopf tensor product
of the representation over $\A(S^4_q)$ with the representation $V_{1/2}$.
From $\A(S^4_q)\simeq\bigoplus_{l\in\N}V_l$ and from the decomposition
$V_l\otimes V_{1/2}\simeq V_{l-\frac{1}{2}}\oplus V_{l+\frac{1}{2}}$
for all $l\in\{1,2,3,\ldots\}$, we deduce that
$
\A(S^4_q)^4\simeq\bigoplus\nolimits_{l\in\N+\frac{1}{2}}(V_l\oplus V_l)
$
and then,
$$
\pi(\A(S^4_q)^4)\simeq\bigoplus\nolimits_{l\in\N+\frac{1}{2}}m^+_lV_l\;,\qquad
\pi^\perp(\A(S^4_q)^4)\simeq\bigoplus\nolimits_{l\in\N+\frac{1}{2}}m^-_lV_l\;,
$$
with multiplicities $m^\pm_l$ to be determined, such that $m_l^++m^-_l=2$.
For $l\in\N+\frac{1}{2}$, the vectors 
$$
v_l^\pm:=x_2^{l-\frac{1}{2}}(1\pm x_0,\pm q^3x_2,\mp qx_1,0) \;.
$$
are highest weight vectors, being annihilated by both $E_1$ and $E_2$, and
have weight $(\frac{1}{2},l)$. Furthermore, $v^+_l(1-e)=v^-_le=0$. Thus,
$v^+_l\in\pi(\A(S^4_q)^4)$ and $v^-_l\in\pi^\perp(\A(S^4_q)^4)$.

Then in both modules $\pi(\A(S^4_q)^4)$ and $\pi^\perp(\A(S^4_q)^4)$
each representation $V_l$, $l\in\N+\frac{1}{2}$, appears with multiplicity
$m_l^\pm\geq 1$. Since $m_l^++m^-_l=2$, we deduce that $m_l^\pm=1$ for all
$l\in\N+\frac{1}{2}$.
\end{prova}

\section{Equivariant representations of $\A(S^4_q)$}\label{sec:eqreps}
Next, we construct $\Uq$-equivariant representations of $\A(S^4_q)$ which
classically correspond to the left regular and chiral spinor representations.
The representation spaces will be (the closure of) $\bigoplus_{l\in\N}V_l$
and $\bigoplus_{l\in\N+\frac{1}{2}}V_l$.

Equivariance of a representation means that it is a representation of the crossed
product algebra $\A(S^4_q)\rtimes\Uq$. The latter is defined by the crossed relations
$ha=(h_{(1)}\az a)h_{(2)}$ for all $a\in\A(S^4_q)$ and $h\in \Uq$; explicitly, the
relations between the generators read:

\begin{equation}\label{eq:crossA}\begin{array}{ccccc}
[K_1,x_0]=0\;, && K_1x_1=qx_1K_1\;, && K_1x_2=x_2K_1\;, \\
\rule{0pt}{18pt}
[K_2,x_0]=0\;, && K_2x_1=q^{-1}x_1K_2\;, && K_2x_2=qx_2K_2\;, \\
\rule{0pt}{18pt}
[E_1,x_0]=q^{-1/2}x_1K_1\;, && E_1x_1=q^{-1}x_1E_1\;, && E_1x_2=x_2E_1\;, \\
\rule{0pt}{18pt}
[F_1,x_0]=-q^{-1/2}K_1x_1^*\;, && F_1x_1=q^{-1}x_1F_1+q^{1/2}[2]x_0K_1\;, && F_1x_2=x_2F_1\;, \\
\rule{0pt}{18pt}
[E_2,x_0]=0\;, && E_2x_1=qx_1E_2+x_2K_2\;, && E_2x_2=q^{-1}x_2E_2\;, \\
\rule{0pt}{18pt}
[F_2,x_0]=0\;, && F_2x_1=qx_1F_2\;, && F_2x_2=q^{-1}x_2F_2+x_1K_2\;.\!\!\!\!
\end{array}\end{equation}
In the previous section we proved that on $\bigoplus_{l\in\N}V_l$ there is at
least one equivariant representation, the left regular one, and that on
$\bigoplus_{l\in\N+\frac{1}{2}}V_l$ there are at least two equivariant
representations, corresponding to the projective modules $\A(S^4_q)^4e$
and $\A(S^4_q)^4(1-e)$. In this section we'll prove that on such spaces
there are no other equivariant representations besides the ones just
mentioned.

Let us denote with $\ket{l,m_1,m_2;j}$ the basis of the representation
space $V_l$ of $\Uq$ as discussed in Sect.~\ref{sec:rep}.
From the first two lines of (\ref{eq:crossA}) we deduce that 
\begin{subequations}\label{eq:coeff}
\begin{align}
x_0\ket{l,m_1,m_2;j} &=\sum_{l',j'}A_{j,j',l,l'}^{m_1,m_2}\ket{l',m_1,m_2;j'} \;, \\
x_1\ket{l,m_1,m_2;j} &=\sum_{l',j'}B_{j,j',l,l'}^{m_1,m_2}\ket{l',m_1+1,m_2;j'} \;, \\
x_2\ket{l,m_1,m_2;j} &=\sum_{l',j'}C_{j,j',l,l'}^{m_1,m_2}\ket{l',m_1,m_2+1;j'} \;,\label{eq:coeffC}
\end{align}
\end{subequations}
with coefficients to be determined. Notice that from the crossed relations
\begin{align*}
x_1\ket{l,m_1,m_2;j} &=(F_2x_2-q^{-1}x_2F_2)K_2^{-1}\ket{l,m_1,m_2;j} \;,\\
x_0\ket{l,m_1,m_2;j} &=q^{-1/2}[2]^{-1}(F_1x_1-q^{-1}x_1F_1)K_1^{-1}\ket{l,m_1,m_2;j} \;,
\end{align*}
the matrix coefficients of $x_0$ and $x_1$ can be expressed in terms
of the coefficients of $x_2$.

\begin{lemma}
Let $k\in\N$. The following formul{\ae} hold:
\begin{subequations}\label{eq:nul}
\begin{align}
F_1^k\ket{l,m_1,m_2;j} &=
  \left\{\begin{array}{ll}
  =0 &\quad\mathrm{if}\;k>j+m_1 \\
  \neq 0&\quad\mathrm{if}\;k\leq j+m_1
  \end{array}\right. \;\; , \label{eq:nulA} \\
E_1^k\ket{l,m_1,m_2;j} &=
  \left\{\begin{array}{ll}
  =0 &\quad\mathrm{if}\;k>j-m_1 \\
  \neq 0&\quad\mathrm{if}\;k\leq j-m_1
  \end{array}\right. \;\; . \label{eq:nulB}
\end{align} 
\end{subequations} 
\end{lemma}

\begin{prova}
By direct computation:
\begin{align*}
F_1^k\ket{l,m_1,m_2;j} &=\sqrt{[j+m_1][j+m_1-1]\ldots [j+m_1-k+1]}\,\times \\ 
& \qquad \times\sqrt{[j-m_1+1][j-m_1+2]\ldots [j-m_1+k]}\,\ket{l,m_1-k,m_2;j}\;.
\end{align*}
The second square root is always different from zero since the
$q$-analogues are in increasing order and $j-m_1+1\geq 1$. In
the first square root $q$-analogues are in decreasing order and
are all different from zero if and only if $j+m_1-k+1\geq 1$.
This proves Eq.~(\ref{eq:nulA}).

In the same way one establishes (\ref{eq:nulB}) by computing that
\begin{align*}
E_1^k\ket{l,m_1,m_2;j} &=\sqrt{[j-m_1][j-m_1-1]\ldots [j-m_1-k+1]}\,\times \\ 
& \qquad
\times\sqrt{[j+m_1+1][j+m_1+2]\ldots [j+m_1+k]}\,\ket{l,m_1+k,m_2;j}\;.
\end{align*}
\end{prova}

\begin{lemma}\label{lemma:A}
The coefficients in (\ref{eq:coeff}) satisfy:
$$
A^{m_1,m_2}_{j,j',l,l'}=B^{m_1,m_2}_{j,j',l,l'}=0\;\;\mathrm{if}\;\;|j-j'|>1\;,
\qquad C^{m_1,m_2}_{j,j',l,l'}=0\;\;\mathrm{if}\;\;j'\neq j\;.
$$
\end{lemma}

\begin{prova}
From (\ref{eq:crossA}), (\ref{eq:nulA}) and (\ref{eq:nulB}) we derive:
\begin{align*}
E_1^{j-m_1+1}x_1\ket{l,m_1,m_2;j} &=q^{-j+m_1-1}x_1E_1^{j-m_1+1}\ket{l,m_1,m_2;j}=0\;, \\
F_1^{j'+m_1+2}x_1^*\ket{l',m_1+1,m_2;j'} &=q^{j'+m_1+2}x_1^*F_1^{j'+m_1+2}\ket{l',m_1+1,m_2;j'}=0\;.
\end{align*}
We expand the left hand sides and use the independence of the vectors
$E_1^{j-m_1+1}\ket{l',m_1+1,m_2;j'}$ and
$F_1^{j'+m_1+2}\ket{l,m_1,m_2;j}$ to arrive at the conditions:
\begin{align*}
B_{j,j',l,l'}^{m_1,m_2}\Big\{E_1^{j-m_1+1}\ket{l',m_1+1,m_2;j'}\Big\}  &=0\;, \\
\bar{B}_{j,j',l,l'}^{m_1,m_2}\Big\{F_1^{j'+m_1+2}\ket{l,m_1,m_2;j}\Big\} &=0\;.
\end{align*}
By (\ref{eq:nulB}) the graph parenthesis in the first line is different from zero
if $j-m_1+1\leq j'-m_1-1$, i.e.~$B_{j,j',l,l'}^{m_1,m_2}$ must be zero
if $j'\geq j+2$. By (\ref{eq:nulA}) the graph parenthesis in the second line is different
from zero if $j'+m_1+2\leq j+m_1$, i.e.~$\bar{B}_{j,j',l,l'}^{m_1,m_2}$ must be zero
if $j'\leq j-2$. This proves $1/3$ of the statement
\begin{equation*}
B_{j,j',l,l'}^{m_1,m_2}=0\quad\forall\;j'\notin\{j-1,j,j+1\}\;.
\end{equation*}
A similar argument applies to $x_0$. From the coproduct of $E_1^n$ we deduce:
\begin{equation*}
E_1^nx_0 =\sum_{k=0}^n\sqbn{n}{k}(E_1^k\az x_0)E_1^{n-k}K_1^k
       =x_0E_1^n-[n]q^{-1/2}x_1E_1^{n-1}K_1\;.
\end{equation*}
This implies that $E_1^{j-m_1+2}x_0\ket{l,m_1,m_2;j}=0$ and
$F_1^{j'+m_1+2}x_0\ket{l',m_1,m_2;j}=0$. From these conditions
we deduce that also $x_0$ shift $j$ by $\{0,\pm 1\}$ only.

Finally, let $\mathcal{C}_1$ be the Casimir element in Eq.~(\ref{eq:Cdef}). Then
$[\mathcal{C}_1,x_2]=0$ and from (\ref{eq:Caz}) we deduce that $x_2$ is diagonal
on the index $j$.
\end{prova}

\begin{lemma}\label{le:B}
The coefficients in (\ref{eq:coeffC}) satisfy
$$
C^{m_1,m_2}_{j,j',l,l'}=0\;\;\mathrm{if}\;\;|l-l'|>1\;\;\textup{or if}\;\;
|l-l'|=0\;\;\mathrm{and}\;\;l\in\N\;.
$$
\end{lemma}
\begin{prova}
The elements $\{x_i,x_i^*\}$ are a basis of the irreducible representation $V_1$.
Covariance of the action tells that $x_i\ket{l,m_1,m_2;j}$ and
$x_i^*\ket{l,m_1,m_2;j}$ are a basis of the tensor representation
$V_1\otimes V_l$. Equations (14--15) in Chapter 7 of~\cite{Kli} tell that
$V_1\otimes V_l\simeq V_{l-1}\oplus V_{l+1}$ if $l\in\N$
and that $V_1\otimes V_l\simeq V_{l-1}\oplus V_l\oplus V_{l+1}$ if $l\in\N+\frac{1}{2}$
(with $V_{l-1}$ omitted if $l-1<0$). This Clebsh-Gordan decomposition
tells that $x_2\ket{l,m_1,m_2;j}$ is in the linear span of the basis
vectors $\ket{l',m_1,m_2+1;j}$ with $l'-l=\pm 1$ if $l\in\N$ or
with $l'-l=0,\pm 1$ if $l\in\N+\frac{1}{2}$. This concludes the proof of the lemma.
\end{prova}

\subsection{Computing the coefficients of $x_2$}
By Lemma~\ref{le:B}, we have to consider only the cases $j'=j$, $|l'-l|\leq 1$ if $l\in\N+\frac{1}{2}$
or $|l'-l|=1$ if $l\in\N$.
The condition $[E_1,x_2]=0$ implies that $C_{j,j,l,l'}^{m_1,m_2}=:C_{j,l,l'}^{m_2}$ is
independent on $m_1$.
Equations $(E_2x_2-q^{-1}x_2E_2)\ket{l,-j,m_2;j}=0$ and
$(F_2x_2^*-qx_2^*F_2)\ket{l',j,m_2+1;j}=0$ imply, respectively:{\small
\begin{gather*}
C^{m_2}_{j,l,l'}\sqrt{[l'-j-m_2-1+\epsilon'][l'+j+m_2+4+\epsilon']}\,=
C^{m_2+1}_{j+1,l,l'}q^{-1}\sqrt{[l-j-m_2+\epsilon][l+j+m_2+3+\epsilon]}\;, \\
C^{m_2}_{j,l,l'}\sqrt{[l+j-m_2+3-\epsilon][l-j+m_2-\epsilon]}\,=
C^{m_2-1}_{j+1,l,l'}q\sqrt{[l'+j-m_2+2-\epsilon'][l'-j+m_2+1-\epsilon']}\;,
\end{gather*}
}with $\epsilon,\epsilon'\in\{0,\pm\frac{1}{2}\}$ determined by the conditions
$l+\epsilon-j-m_2\in 2\N$ and $l'-\epsilon'-j-m_2\in 2\N$.
Notice that if $l'-l\in 2\N+1$ then $\epsilon'=\epsilon$, while if
$l'-l\in 2\N$ then $\epsilon'=-\epsilon$.
Looking at the cases $l'-l=\pm 1$, we deduce that
$$
\frac{ q^{-\frac{1}{2}(j+m_2)} }{\sqrt{[l+j+m_2+3+\epsilon]}}\, C_{j,l,l+1}^{m_2}
\qquad\mathrm{and}\qquad
\frac{ q^{-\frac{1}{2}(j+m_2)} }{\sqrt{[l-j-m_2+\epsilon]}}\, C_{j,l,l-1}^{m_2} 
$$
depend on $j+m_2$ only through their parity (i.e.~they depend only on the value of $\epsilon$). Similarly,
$$
\frac{ q^{\frac{1}{2}(j-m_2)} }{\sqrt{[l-j+m_2+2-\epsilon]}}\, C_{j,l,l+1}^{m_2} 
\qquad\mathrm{and}\qquad
\frac{ q^{\frac{1}{2}(j-m_2)} }{\sqrt{[l+j-m_2+1-\epsilon]}}\, C_{j,l,l-1}^{m_2} 
$$
depend on $j-m_2$ only through their parity.
Combining these informations, we deduce that the following elements do not depend on
the exact value of $j$, $m_2$, but only on the value of $\epsilon$,
\begin{align*}
\frac{ q^{-m_2} }{\sqrt{[l+j+m_2+3+\epsilon][l-j+m_2+2-\epsilon]}}\, C_{j,l,l+1}^{m_2} &=:C_{l,l+1}(\epsilon) \;, \\
\frac{ q^{-m_2} }{\sqrt{[l-j-m_2+\epsilon][l+j-m_2+1-\epsilon]}}\, C_{j,l,l-1}^{m_2} 
 &=:C_{l,l-1}(\epsilon) \;.
\end{align*}
If $l\in\N$ there are no other coefficients $C_{j,l,l'}^{m_2}$ to compute.
If $l\notin\N$, we have to compute also $C_{j,l,l}^{m_2}$.
In this case $\epsilon'=-\epsilon$ and we get:{\small
\begin{gather*}
C^{m_2}_{j,l,l}\sqrt{[l-j-m_2-1-\epsilon][l+j+m_2+4-\epsilon]}\,=
C^{m_2+1}_{j+1,l,l}q^{-1}\sqrt{[l-j-m_2+\epsilon][l+j+m_2+3+\epsilon]}\;, \\
C^{m_2}_{j,l,l}\sqrt{[l+j-m_2+3-\epsilon][l-j+m_2-\epsilon]}\,=
C^{m_2-1}_{j+1,l,l}q\sqrt{[l+j-m_2+2+\epsilon][l-j+m_2+1+\epsilon]}\;.
\end{gather*}
}Again, looking at the two cases $\epsilon=\pm\frac{1}{2}$ we deduce that
$$
\frac{ q^{-\frac{1}{2}(j+m_2)} }{\sqrt{[l+\frac{1}{2}-j-m_2]}}\, C_{j,l,l}^{m_2} 
\;\;\;\mathrm{if}\;\epsilon=\tfrac{1}{2}
\quad\mathrm{and}\quad
\frac{ q^{-\frac{1}{2}(j+m_2)} }{\sqrt{[l+\frac{1}{2}+j+m_2+2]}}\, C_{j,l,l}^{m_2} 
\;\;\;\mathrm{if}\;\epsilon=-\tfrac{1}{2}
$$
do not depend on $j+m_2$ (this time $\epsilon$ is fixed, so the parity of $j+m_2$ is fixed). Similarly,
$$
\frac{ q^{\frac{1}{2}(j-m_2)} }{\sqrt{[l+\frac{1}{2}-j+m_2+1]}}\, C_{j,l,l}^{m_2} 
\;\;\;\mathrm{if}\;\epsilon=\tfrac{1}{2}
\quad\mathrm{and}\quad
\frac{ q^{\frac{1}{2}(j-m_2)} }{\sqrt{[l+\frac{1}{2}+j-m_2+1]}}\, C_{j,l,l}^{m_2} 
\;\;\;\mathrm{if}\;\epsilon=-\tfrac{1}{2}
$$
do not depend on $j-m_2$.
Combining these informations, we deduce that the following element does not depend on
the exact value of $j$, $m_2$, but only on the value of $\epsilon$:
$$
\frac{ q^{-m_2} }{\sqrt{[l-2\epsilon j-m_2+1-\epsilon]
[l-2\epsilon j+m_2+2-\epsilon]}}\, C_{j,l,l}^{m_2} =:C_{l,l}(\epsilon) \;.
$$
The denominator of the left-hand side is just $[2j][2j+2]b_l(j,m_2)$
with $b_l$ the coefficient in Eq.~(\ref{eq:abc}).
The formula $C_{j,l,l}^{m_2}=q^{m_2}[2j][2j+2]b_l(j,m_2)C_{l,l}(\epsilon)$ is
valid for all $l$, since $b_l(j,m_2)$ vanish for $l$ integer.

Summarizing, we find that
\begin{subequations}\label{eq:Cjllm}
\begin{align}
C_{j,l,l+1}^{m_2} &=q^{m_2}\sqrt{[l+j+m_2+3+\epsilon][l-j+m_2+2-\epsilon]} \, C_{l,l+1}(\epsilon) \;, \\
C_{j,l,l}^{m_2} &=q^{m_2}[2j][2j+2]b_l(j,m_2) \, C_{l,l}(\epsilon) \;, \\
C_{j,l,l-1}^{m_2} &=q^{m_2}\sqrt{[l-j-m_2+\epsilon][l+j-m_2+1-\epsilon]} \, 
C_{l,l-1}(\epsilon)\;,
\end{align}
\end{subequations}
with coefficients $C_{l,l'}(\epsilon)$ to be determined.

\subsection{Computing the coefficients of $x_1$}
From Lemma~\ref{lemma:A}, we have to consider only the three cases 
$j'=j,j\pm 1$.
Using equation $E_1x_1=q^{-1}x_1E_1$ we get, 
\begin{align*}
\frac{ q^{-m_1} }{\sqrt{[j+m_1+1][j+m_1+2]}}\, B_{j,j+1,l,l'}^{m_1,m_2}  &=
\frac{ q^{-m_1-1} }{\sqrt{[j+m_1+2][j+m_1+3]}}\, B_{j,j+1,l,l'}^{m_1+1,m_2} \;, \\
\frac{ q^{-m_1} }{\sqrt{[j-m_1][j+m_1+1]}}\, B_{j,j,l,l'}^{m_1,m_2}  &=
\frac{ q^{-m_1-1} }{\sqrt{[j-m_1-1][j+m_1+2]}}\, B_{j,j,l,l'}^{m_1+1,m_2} \;, \\
\frac{ q^{-m_1} }{\sqrt{[j-m_1][j-m_1-1]}}\, B_{j,j-1,l,l'}^{m_1,m_2} &=
\frac{ q^{-m_1-1} }{\sqrt{[j-m_1-1][j-m_1-2]}}\, B_{j,j-1,l,l'}^{m_1+1,m_2} \;.
\end{align*}
We see that the left hand sides of these three equations are independent of $m_1$, and call:
\begin{subequations}\label{eq:Bm}
\begin{align}
B_{j,j+1,l,l'}^{m_1,m_2} &=:q^{m_1}\sqrt{[j+m_1+1][j+m_1+2]}\,B_{j,j+1,l,l'}^{m_2} \;, \\
B_{j,j,l,l'}^{m_1,m_2}    &=:q^{m_1}\sqrt{[j-m_1][j+m_1+1]}\,B_{j,j,l,l'}^{m_2} \;, \\
B_{j,j-1,l,l'}^{m_1,m_2} &=:q^{m_1}\sqrt{[j-m_1][j-m_1-1]}\,B_{j,j-1,l,l'}^{m_2} \;.
\end{align}
\end{subequations}
Imposing the condition $x_1K_2=F_2x_2-q^{-1}x_2F_2$ on the subspace
spanned by $\ket{l,j,m_2;j}$ (so $m_1=j$ and $B_{j,j,l,l'}^{m_1,m_2}=B_{j,j-1,l,l'}^{m_1,m_2}=0$
on this subspace) we get:
$$
q^{m_2}B^{m_2}_{j,j+1,l,l'}=c_{l'}(j+1,m_2)C^{m_2}_{j,l,l'}-q^{-1}c_l(j+1,m_2-1)C^{m_2-1}_{j+1,l,l'}\;.
$$
From this we deduce that, that since coefficients $C_{j,l,l'}^{m_2}$ vanish
for $|l-l'|>1$, also $B_{j,j+1,l,l'}^{m_2}$ is zero in these cases.
In the remaining three cases $l'=l,l\pm 1$, using equation
(\ref{eq:Cjllm}) we get:
\begin{subequations}\label{eq:Bjp}
\begin{align}
B_{j,j+1,l,l+1}^{m_2} &=(-1)^{2\epsilon} q^{l-j+m_2-\epsilon}
  \sqrt{\frac{[l+j+m_2+3+\epsilon][l+j-m_2+3-\epsilon]}{[2(j+|\epsilon|)+1]
  [2(j-|\epsilon|)+3]}}\, C_{l,l+1}(\epsilon) \;, \\
B_{j,j+1,l,l}^{m_2} &=-2\epsilon q^{2\epsilon\,l-j+m_2-2+3\epsilon}
\frac{\sqrt{[l+\frac{1}{2}+j-2\epsilon\,m_2+2][l+\frac{1}{2}-j-2\epsilon\,m_2]}}{[2j+2]}
  \, C_{l,l}(\epsilon) \;, \\
B_{j,j+1,l,l-1}^{m_2} &=(-1)^{2\epsilon+1} q^{-l-j+m_2-3+\epsilon}
  \sqrt{\frac{[l-j+m_2-\epsilon][l-j-m_2+\epsilon]}{[2(j+|\epsilon|)+1]
 [2(j-|\epsilon|)+3]}}\, C_{l,l-1}(\epsilon) \;.
\end{align}
\end{subequations}
Imposing $qx_1^*K_2=x_2^*E_2-q^{-1}E_2x_2^*$ on the subspace spanned by
$\ket{l',-j+1,m_2;j-1}$  (so $B_{j,j,l,l'}^{m_1,m_2}=B_{j,j+1,l,l'}^{m_1,m_2}=0$
on this subspace) we get:
$$
q^{m_2}B^{m_2}_{j,j-1,l,l'}=a_{l'}(j-1,m_2)C^{m_2}_{j,l,l'}-q^{-1}a_l(j-1,m_2-1)C^{m_2-1}_{j-1,l,l'}\;.
$$
We deduce that $B_{j,j-1,l,l'}^{m_2}$ vanishes if $|l-l'|>1$, while in the three remaining
cases $l'=l,l\pm 1$ using Eq.~(\ref{eq:Cjllm}) we get:
\begin{subequations}\label{eq:Bjm}
\begin{align}
& B_{j,j-1,l,l+1}^{m_2}=q^{l+j+m_2+1+\epsilon}
  \sqrt{\frac{[l-j-m_2+2+\epsilon][l-j+m_2+2-\epsilon]}{[2(j+|\epsilon|)-1][2(j-|\epsilon|)+1]}}\, C_{l,l+1}(\epsilon) \;, \\
& B_{j,j-1,l,l}^{m_2}=-2\epsilon q^{-2\epsilon\,l+j+m_2-1-3\epsilon}
\frac{\sqrt{[l+\frac{1}{2}+j+2\epsilon\,m_2+1][l+\frac{1}{2}-j+2\epsilon\,m_2+1]}}{[2j]}\, C_{l,l}(\epsilon)  \;, \\
& B_{j,j-1,l,l-1}^{m_2}=-q^{-l+j+m_2-2-\epsilon}
  \sqrt{\frac{[l+j+m_2+1+\epsilon][l+j-m_2+1-\epsilon]}{[2(j+|\epsilon|)-1][2(j-|\epsilon|)+1]}}\, C_{l,l-1}(\epsilon)\;.
\end{align}
\end{subequations}
Moreover, the condition 
$\inner{l',j,m_2;j|x_1K_2+q^{-1}x_2F_2-F_2x_2|l,j-1,m_2;j}=0$
implies that
\begin{equation}\label{eq:xxx}
q^{m_2}B^{m_2}_{j,j,l,l'}=b_{l'}(j,m_2)C^{m_2}_{j,l,l'}-q^{-1}b_l(j,m_2-1)C^{m_2-1}_{j,l,l'}\;.
\end{equation}
A further elaboration on these coefficients is postponed to after the following section.

\subsection{Computing the coefficients of $x_0$}
The condition $q^{1/2}[2]x_0K_1=F_1x_1-q^{-1}x_1F_1$ implies:
$$
q^{m_1+\frac{1}{2}}[2]A_{j,j',l,l'}^{m_1,m_2}=\sqrt{[j'-m_1][j'+m_1+1]}\, B_{j,j',l,l'}^{m_1,m_2}-q^{-1}\sqrt{[j+m_1][j-m_1+1]}\, B_{j,j',l,l'}^{m_1-1,m_2}\;.
$$
In the three non-trivial cases $j'-j=1,0,-1$, using (\ref{eq:Bm}), we get:
\begin{subequations}\label{eq:Am}
\begin{align}
A_{j,j+1,l,l'}^{m_1,m_2} &=q^{j+m_1-\frac{1}{2}}\sqrt{[j+m_1+1][j-m_1+1]}\,B^{m_2}_{j,j+1,l,l'} \;, \\
A_{j,j,l,l'}^{m_1,m_2} &=[2]^{-1}q^{-2-\frac{1}{2}}\big(q^{j+m_1+1}[2][j-m_1]-[2j]\big)B_{j,j,l,l'}^{m_2}\;, \\
A_{j,j-1,l,l'}^{m_1,m_2} &=-q^{-j+m_1-1-\frac{1}{2}}\sqrt{[j+m_1][j-m_1]}\,B^{m_2}_{j,j-1,l,l'} \;.
\end{align}
\end{subequations}
The hermiticity condition $x_0=x_0^*$ means that $A^{m_1,m_2}_{j,j+1,l,l'}=\bar{A}^{m_1,m_2}_{j+1,j,l'l}$
and $A^{m_1,m_2}_{j,j,l,l'}=\bar{A}^{m_1,m_2}_{j,j,l'l}$. Thus, from
(\ref{eq:Am}) it follows that:
$$
B^{m_2}_{j+1,j,l'l}=-q^{2j+2}\bar{B}^{m_2}_{j,j+1,l,l'} \;,\qquad
B^{m_2}_{j,j,l,l'}=\bar{B}^{m_2}_{j,j,l'l} \;.
$$
Using (\ref{eq:Bjp}), the first equation turns out to be equivalent to
the following conditions:
\begin{subequations}\label{eq:sol}
\begin{equation}
C_{l+1,l}(\epsilon)=(-1)^{2\epsilon}q^{2l+4}\bar{C}_{l,l+1}(\epsilon)\;,\qquad
C_{l,l}(\epsilon)=\bar{C}_{l,l}(-\epsilon)\;.
\end{equation}
The second of equation together with (\ref{eq:xxx}) implies:
$$
b_{l'}(j,m_2)C^{m_2}_{j,l,l'}-q^{-1}b_l(j,m_2-1)C^{m_2-1}_{j,l,l'}=
b_l(j,m_2)C^{m_2}_{j,l'l}-q^{-1}b_{l'}(j,m_2-1)C^{m_2-1}_{j,l'l}\;.
$$
That is, using (\ref{eq:Cjllm}):
\begin{equation}
C_{l,l+1}(\epsilon)=C_{l,l+1}(-\epsilon)\;,\qquad
C_{l,l}(\epsilon)=\bar{C}_{l,l}(\epsilon)\;.
\end{equation}
\end{subequations}

\subsection{Again the coefficients of $x_1$}
Now, using (\ref{eq:sol}) together with (\ref{eq:xxx}) we are able to compute
the last coefficients. Notice that from (\ref{eq:abc}) the coefficients $b_l$ vanish
if $\epsilon=0$ (i.e.~in the left regular representation), and then from (\ref{eq:xxx})
$B^{m_2}_{j,j,l,l'}$ vanish too if $\epsilon=0$.
Moreover, from Lemma \ref{lemma:A} $B^{m_2}_{j,j,l,l'}$ vanish also if $|l-l'|>1$.
In the three cases $l'=l,l\pm 1$, using Eq.~(\ref{eq:Cjllm}) we get:
\begin{align*}
B_{j,j,l,l+1}^{m_2} &=2|\epsilon|[2]q^{l+m_2+1+\epsilon(2j+1)}
  \frac{\sqrt{[l+2\epsilon\,j-m_2+2+\epsilon][l-2\epsilon\,j+m_2+2-\epsilon]}}{[2j][2j+2]}\, C_{l,l+1}(\epsilon) \;, \\
B_{j,j,l,l}^{m_2} &=\frac{2|\epsilon|}{[2j][2j+2]}\Big\{
                  [l-\epsilon(2j+1)-m_2+1][l-\epsilon(2j+1)+m_2+2]+ \\
&\phantom{=\frac{2|\epsilon|\,C_{l,l}(\epsilon)}{[2j][2j+2]}\, C_{l,l}(\epsilon)\Big\{}
                 -q^{-2}[l+\epsilon(2j+1)-m_2+2][l+\epsilon(2j+1)+m_2+1]\Big\} \\
         &=-\frac{2|\epsilon|}{[2j][2j+2]}\,
           \frac{q^{-2\epsilon(2j+1)}[2l+4]-q^{2\epsilon(2j+1)}[2l+2]-[2]q^{2m_2}}{1-q^2}\, C_{l,l}(\epsilon) \;, \\
B_{j,j,l,l-1}^{m_2} &=-2|\epsilon|[2]q^{-l+m_2-2-\epsilon(2j+1)}
  \frac{\sqrt{[l+2\epsilon\,j+m_2+1+\epsilon][l-2\epsilon\,j-m_2+1-\epsilon]}}{[2j][2j+2]}\, C_{l,l-1}(\epsilon) \;.
\end{align*}
We have inserted the factor $2|\epsilon|$, so that the expressions remain valid also
when $\epsilon=0$.

\subsection{The condition on the radius}
Orbits for $SO(5)$ are spheres of arbitrary radius, equivariance
alone  not imposing constraints on the radius. Similarly, for the quantum spheres 
one has to impose a constraint on the radius to determine the coefficients of the representation.
In fact, this will determine $C_{l,l+1}(0)$, $C_{l,l+1}(\frac{1}{2})$ and
$C_{l,l}(\frac{1}{2})$ only up to a phase.
Different choices of the phases correspond to unitary equivalent
representations and without losing generality we choose
$C_{l,l'}(\epsilon)\in\R$. A possible expression for the radius is 
$q^8x_0^2+q^4x_1^*x_1+x_2^*x_2$ which we constrain to be equal to 1. 
Let then, 
$$
r(l,m_1,m_2;j):=\inner{l,m_1,m_2;j\big|\big(q^8x_0^2+q^4x_1^*x_1+x_2^*x_2
\big)\big|l,m_1,m_2;j}\;.
$$
All these matrix coefficients must be $1$. In particular,
for $l\in\N$ the condition $\,r(l,0,l;0)=1\,$ implies (up to a phase) that
\begin{equation}\label{eq:lint}
C_{l,l+1}(0)=\frac{q^{-l-3/2}}{\sqrt{[2l+3][2l+5]}}\;.
\end{equation}
For $l\in\N+\frac{1}{2}$ we first require that
$\,r(l,\frac{1}{2},l;\frac{1}{2})=r(l,-\frac{1}{2},l;\frac{1}{2})\,$
obtaining two possibilities:
$$
C_{l,l}(\tfrac{1}{2})=\pm\frac{[2]q^{l+2}}{[2l+2]}C_{l,l+1}(\tfrac{1}{2})\;.
$$
Then imposing $\,r(l,\tfrac{1}{2},l;\tfrac{1}{2})=1\,$, yields (up to a phase)
\begin{equation}
C_{l,l+1}(\tfrac{1}{2})=\frac{q^{-l-3/2}}{[2l+4]}\;, 
\end{equation}
hence,
\begin{equation}\label{eq:twocases}
C_{l,l}(\tfrac{1}{2})=\pm\frac{q^{1/2}[2]}{[2l+2][2l+4]}\;.
\end{equation}
With these, all the coefficients are completely determined.

\subsection{Explicit form of the representations}
Let us recall what we know on the equivariant representations of
the algebra $\A(S^4_q)$.

By the decomposition $\A(S^4_q)\simeq\bigoplus_{l\in\N}V_l$ into
irreducible representations of $U_q(so(5))$, there exists (at least)
one representation of $\A(S^4_q)\rtimes U_q(so(5))$ on the vector
space $\bigoplus_{l\in\N}V_l$ extending the representation
$\bigoplus_{l\in\N}\sigma_l$ of $U_q(so(5))$.
As we computed above, the equivariance uniquely
determines (for $l\in\N$, up to unitary equivalence) the matrix
coefficients of the representation, whose expression is
characterized by (\ref{eq:lint}).
On the other hand, by Proposition \ref{prop:main} there are
(at least) two inequivalent representations of $\A(S^4_q)\rtimes U_q(so(5))$ on the vector
space $\bigoplus_{l\in\N+\frac{1}{2}}V_l$ extending the representation
$\bigoplus_{l\in\N+\frac{1}{2}}\sigma_l$ of $U_q(so(5))$. These correspond,
by Lemma \ref{lemma:ptwo}, to the projective modules $\A(S^4_q)^4e$ and
$\A(S^4_q)^4(1-e)$, with $e$ the idempotent in Eq.~(\ref{eq:P}).
The computation above (for $l\in\N+\frac{1}{2}$), which culminates
in Eq.~(\ref{eq:twocases}), tells us that there are only
two possibilities for the matrix coefficients (up to unitary equivalence).
Therefore, the two possible choices in (\ref{eq:twocases}) must correspond
to the inequivalent representations associated with the projective
modules $\A(S^4_q)^4e$ and $\A(S^4_q)^4(1-e)$.

Let us summarize these results in the following
two theorems, which correspond to the scalar (i.e.~left regular)
and chiral spinor representations, respectively.

\begin{thm}
The vector space
$\A(S^4_q)$ has orthonormal basis $\ket{l,m_1,m_2;j}$ with,
$$
l\in\N\;, \qquad j=0,1,\ldots,l\;,\qquad
j-|m_1|\in\N\;,\qquad
l-j-|m_2|\in 2\N\;.
$$
We call $L^2(S^4_q)$ the Hilbert space completion of $\A(S^4_q)$. Modulo a
unitary equivalence, the left regular representation is given by{\allowdisplaybreaks[0]
\begin{align*}
x_0\ket{l,m_1,m_2;j} &=A_{j,m_1}C^+_{l,j,m_2}\ket{l+1,m_1,m_2;j+1} \\
                     &+A_{j,m_1}C^-_{l,j,m_2}\ket{l-1,m_1,m_2;j+1} \\
                     &+A_{j-1,m_1}C^-_{l+1,j-1,m_2}\ket{l+1,m_1,m_2;j-1} \\
                     &+A_{j-1,m_1}C^+_{l-1,j-1,m_2}\ket{l-1,m_1,m_2;j-1} \;, \\
x_1\ket{l,m_1,m_2;j} &=B^+_{j,m_1}C^+_{l,j,m_2}\ket{l+1,m_1+1,m_2;j+1} \\
                     &+B^+_{j,m_1}C^-_{l,j,m_2}\ket{l-1,m_1+1,m_2;j+1} \\
                     &+B^-_{j,m_1}C^-_{l+1,j-1,m_2}\ket{l+1,m_1+1,m_2;j-1} \\
                     &+B^-_{j,m_1}C^+_{l-1,j-1,m_2}\ket{l-1,m_1+1,m_2;j-1} \;, \\
x_2\ket{l,m_1,m_2;j} &=D_{l,j,m_2}^+\ket{l+1,m_1,m_2+1;j} \\
                     &+D_{l,j,m_2}^-\ket{l-1,m_1,m_2+1;j} \;,
\end{align*}
}with coefficients
\begin{align*}
A_{j,m_1} &=q^{m_1-1}\sqrt{\frac{[j+m_1+1][j-m_1+1]}{[2j+1][2j+3]}} \;, \\
B^+_{j,m_1} &=q^{-j+m_1-1/2}\sqrt{\frac{[j+m_1+1][j+m_1+2]}{[2j+1][2j+3]}} \;, \\
B^-_{j,m_1} &=-q^{j+m_1+1/2}\sqrt{\frac{[j-m_1][j-m_1-1]}{[2j-1][2j+1]}} \;.
\end{align*}
and
\begin{align*}
C^+_{l,j,m_2} &=q^{m_2-1}\sqrt{\frac{[l+j+m_2+3][l+j-m_2+3]}{[2l+3][2l+5]}} \;, \\
C^-_{l,j,m_2} &=-q^{m_2-1}\sqrt{\frac{[l-j+m_2][l-j-m_2]}{[2l+1][2l+3]}} \;, \\
D_{l,j,m_2}^+ &=q^{-l+m_2-3/2}\sqrt{\frac{[l+j+m_2+3][l-j+m_2+2]}{[2l+3][2l+5]}} \;, \\
D_{l,j,m_2}^- &=q^{l+m_2+3/2}\sqrt{\frac{[l-j-m_2][l+j-m_2+1]}{[2l+1][2l+3]}} \;.
\end{align*}
\end{thm}

The two chiral spinorial representations (corresponding to the
sign $\pm$ in Eq.~(\ref{eq:twocases})) are described in
the following theorem.

\begin{thm}\label{thm:cr}
Let $\HH_\pm$ be two Hilbert spaces with orthonormal basis $\ket{l,m_1,m_2;j}_\pm$, where
$$
l\in\N+\tfrac{1}{2}\;, \qquad j=\tfrac{1}{2},\tfrac{3}{2},\ldots,l\;,\qquad
j-|m_1|\in\N\;,\qquad
l+\tfrac{1}{2}-j-|m_2|\in\N\;.
$$
Let $\epsilon=\pm\frac{1}{2}$ be defined by $l+\epsilon-j-m_2\in 2\N$. On each space $\HH_\pm$
there is an equivariant $*$-representation of $\A(S^4_q)$ defined by:
\begin{align*}
x_0\ket{l,m_1,m_2;j}_\pm &=A^+_{j,m_1}C^+_{l,j,m_2}\ket{l+1,m_1,m_2;j+1}_\pm \\
                     &\mp A^+_{j,m_1}C^0_{l,j,m_2}\ket{l,m_1,m_2;j+1}_\pm \\
                     &+A^+_{j,m_1}C^-_{l,j,m_2}\ket{l-1,m_1,m_2;j+1}_\pm \\
                     &+A^0_{j,m_1}H^+_{l,j,m_2}\ket{l+1,m_1,m_2;j}_\pm \\
                     &\pm A^0_{j,m_1}H^0_{l,j,m_2}\ket{l,m_1,m_2;j}_\pm \\
                     &+A^0_{j,m_1}H^+_{l-1,j,m_2}\ket{l-1,m_1,m_2;j}_\pm \\
                     &+A^+_{j-1,m_1}C^-_{l+1,j-1,m_2}\ket{l+1,m_1,m_2;j-1}_\pm \\
                     &\mp A^+_{j-1,m_1}C^0_{l,j-1,m_2}\ket{l,m_1,m_2;j-1}_\pm \\
                     &+A^+_{j-1,m_1}C^+_{l-1,j-1,m_2}\ket{l-1,m_1,m_2;j-1}_\pm \;, \\
x_1\ket{l,m_1,m_2;j}_\pm &=B^+_{j,m_1}C^+_{l,j,m_2}\ket{l+1,m_1+1,m_2;j+1}_\pm \\
                     &\mp B^+_{j,m_1}C^0_{l,j,m_2}\ket{l,m_1+1,m_2;j+1}_\pm \\
                     &+B^+_{j,m_1}C^-_{l,j,m_2}\ket{l-1,m_1+1,m_2;j+1}_\pm \\
                     &+B^0_{j,m_1}H^+_{l,j,m_2}\ket{l+1,m_1+1,m_2;j}_\pm \\
                     &\pm B^0_{j,m_1}H^0_{l,j,m_2}\ket{l,m_1+1,m_2;j}_\pm \\
                     &+B^0_{j,m_1}H^+_{l-1,j,m_2}\ket{l-1,m_1+1,m_2;j}_\pm \\
                     &+B^-_{j,m_1}C^-_{l+1,j-1,m_2}\ket{l+1,m_1+1,m_2;j-1}_\pm \\
                     &\mp B^-_{j,m_1}C^0_{l,j-1,m_2}\ket{l,m_1+1,m_2;j-1}_\pm \\
                     &+B^-_{j,m_1}C^+_{l-1,j-1,m_2}\ket{l-1,m_1+1,m_2;j-1}_\pm \;, \\
x_2\ket{l,m_1,m_2;j}_\pm &=D_{l,j,m_2}^+\ket{l+1,m_1,m_2+1;j}_\pm \\
                     &\pm D_{l,j,m_2}^0\ket{l,m_1,m_2+1;j}_\pm \\
                     &+D_{l,j,m_2}^-\ket{l-1,m_1,m_2+1;j}_\pm \;,
\end{align*}
with coefficients
\begin{align*}
A^+_{j,m_1} &=q^{m_1-1}\frac{\sqrt{[j+m_1+1][j-m_1+1]}}{[2j+2]} \;, \\
A^0_{j,m_1} &=q^{-2}\frac{q^{j+m_1+1}[2][j-m_1]-[2j]}{[2j][2j+2]} \;, \\
B^+_{j,m_1} &=q^{-j+m_1-1/2}\frac{\sqrt{[j+m_1+1][j+m_1+2]}}{[2j+2]} \;, \\
B^0_{j,m_1} &=(1+q^2)q^{m_1-1/2}\frac{\sqrt{[j-m_1][j+m_1+1]}}{[2j][2j+2]} \;, \\
B^-_{j,m_1} &=-q^{j+m_1+1/2}\frac{\sqrt{[j-m_1][j-m_1-1]}}{[2j]} \;.
\end{align*}
and
\begin{align*}
C^+_{l,j,m_2} &=-q^{m_2-1-\epsilon}\frac{\sqrt{[l+j+m_2+3+\epsilon][l+j-m_2+3-\epsilon]}}{[2l+4]} \;, \\
C^0_{l,j,m_2} &=[4\epsilon]\,q^{2\epsilon\,l+m_2-1+3\epsilon}
    \frac{\sqrt{[l+\frac{1}{2}+j-2\epsilon\,m_2+2][l+\frac{1}{2}-j-2\epsilon\,m_2]}}
    {[2l+2][2l+4]} \;, \\
C^-_{l,j,m_2} &=-q^{m_2-1+\epsilon}\frac{\sqrt{[l-j+m_2-\epsilon][l-j-m_2+\epsilon]}}{[2l+2]} \;, \\
H^+_{l,j,m_2} &=q^{m_2-1+\epsilon(2j+1)}\frac{\sqrt{[l+2\epsilon j-m_2+2+\epsilon]
    [l-2\epsilon j+m_2+2-\epsilon]}}{[2l+4]} \;, \\
H^0_{l,j,m_2} &=\tfrac{[l-\epsilon(2j+1)-m_2+1][l-\epsilon(2j+1)+m_2+2]
       -q^{-2}[l+\epsilon(2j+1)-m_2+2][l+\epsilon(2j+1)+m_2+1]}{[2l+2][2l+4]} \;, \\
D_{l,j,m_2}^+ &=q^{-l+m_2-3/2}\frac{\sqrt{[l+j+m_2+3+\epsilon][l-j+m_2+2-\epsilon]}}{[2l+4]} \;, \\
D_{l,j,m_2}^0 &=[2]q^{m_2+1/2}\frac{\sqrt{[l-2\epsilon j-m_2+1-\epsilon]
    [l-2\epsilon j+m_2+2-\epsilon]}}{[2l+2][2l+4]} \;, \\
D_{l,j,m_2}^- &=-q^{l+m_2+3/2}\frac{\sqrt{[l-j-m_2+\epsilon][l+j-m_2+1-\epsilon]}}{[2l+2]} \;.
\end{align*}
These two representations are inequivalent and correspond to the
projective modules $\A(S^4_q)^4e$ and $\A(S^4_q)^4(1-e)$, with $e$ the
idempotent in Eq.~(\ref{eq:P}).
\end{thm}

\section{The Dirac operator on the orthogonal quantum  $4$-sphere}\label{sec:D}
We start by constructing a non-trivial Fredholm module on the
orthogonal quantum sphere (with different representations a non-trivial
Fredholm module was already constructed in~\cite{LH}).

\begin{prop}\label{prop:notr}
Consider the representations of $\A(S^4_q)$ on $\HH_\pm$ given in Theorem \ref{thm:cr}.
Then, the datum $(\A(S^4_q),\HH,F,\gamma)$ is a $1$-summable even Fredholm module,
where $\HH:=\HH_+\oplus\HH_-$, $\gamma$ is the natural grading and $F\in\B(\HH)$ is defined by
$$
F\ket{l,m_1,m_2;j}_\pm:=\ket{l,m_1,m_2;j}_\mp\;.
$$
This Fredholm module is non-trivial. In particular,
\begin{equation}\label{eq:pair}
\mathrm{ch}^F\!([e]):=
\tfrac{1}{2}\tr_{\HH\otimes\C^2}(\gamma F[F,P])=1\;,
\end{equation}
with $e$ the idempotent defined by Eq.~(\ref{eq:P}).
\end{prop}
\begin{prova}
That $F=F^*$, $F^2=1$ and $\gamma F+F\gamma=0$ is obvious. Then, it is enough to show
that $[F,x_i]\in\mathcal{L}^1(\HH)$ for $i=0,1,2$. From this and the Leibniz rule
it follows that $[F,a]$ is trace class, and then compact, for all $a\in\A(S^4_q)$.

Now, notice that
\begin{align}
[F,x_0]\ket{l,m_1,m_2;j}_\pm=\;&\mp 2A^+_{j,m_1}C^0_{l,j,m_2} \ket{l,m_1,m_2;j+1}_\mp \nonumber \\
                               &\pm 2A^0_{j,m_1}H^0_{l,j,m_2} \ket{l,m_1,m_2;j}_\mp \label{eq:Fx0} \\
                               &\mp 2A^+_{j-1,m_1}C^0_{l,j-1,m_2} \ket{l,m_1,m_2;j-1}_\mp \nonumber \;, \\
[F,x_1]\ket{l,m_1,m_2;j}_\pm=\;&\mp 2B^+_{j,m_1}C^0_{l,j,m_2} \ket{l,m_1+1,m_2;j+1}_\mp \nonumber \\
                               &\pm 2B^0_{j,m_1}H^0_{l,j,m_2} \ket{l,m_1+1,m_2;j}_\mp \nonumber \\
                               &\mp 2B^-_{j,m_1}C^0_{l,j-1,m_2} \ket{l,m_1+1,m_2;j-1}_\mp \;, \nonumber \\
[F,x_2]\ket{l,m_1,m_2;j}_\pm=\;&\pm 2D^0_{l,j,m_2}\ket{l,m_1,m_2+1;j}_\mp \;. \nonumber
\end{align}
All the coefficients appearing in these equations are bounded by $q^{2l}$.
Thus the commutators are trace class and this concludes the first part of the proof.

To prove non-triviality it is enough to prove \eqref{eq:pair}.
Substituting \eqref{eq:P} into \eqref{eq:pair} yields
$$
\mathrm{ch}^F\!([e])
=\tfrac{(1-q^2)^2}{4}\tr_{\HH}(\gamma F[F,x_0])\;.
$$
and in turn, using Eq.~(\ref{eq:Fx0}),
$$
\mathrm{ch}^F\!([e])
=(1-q^2)^2\sum_{l,j,m_1,m_2}A^0_{j,m_1}H^0_{l,j,m_2} \;.
$$
Summing over $m_1$ from $-j$ to $j$ we obtain that
\begin{align*}
\mathrm{ch}^F\!([e])
        &=q^{-3}(1-q^2)^2\sum_{l,j,m_2}\frac{[l+\epsilon(2j+1)-m_2+2]
        [l+\epsilon(2j+1)+m_2+1]}{[2l+2][2l+4][2j][2j+2]}\times \\ &
        \qquad\qquad \qquad\qquad \qquad\qquad\qquad\qquad\qquad
        \times
        \sum_{m_1}\big\{q^{2j+2}+q^{-2j}-[2]q^{2m_1+1}\big\} \\
        &=\sum_{l,j}\frac{(2j+1)(q^{2j+1}+q^{-2j-1})-[2][2j+1]}
        {[2l+2][2l+4][2j][2j+2]} \times \\ &
        \qquad\qquad \qquad\qquad
        \times\sum_{m_2}\Big\{q^{2l+2\epsilon(2j+1)+3}+q^{-2l-2\epsilon(2j+1)-3}
        -q^{2m_2-1}-q^{-2m_2+1}\Big\}\;.
\end{align*}

The sum over $m_2$ requires additional care.
For $\epsilon$ fixed, $l-\epsilon-j+m_2=0,2,4,\ldots,2(l-j)$. If we call
$2i:=l-\epsilon-j+m_2$ and sum first over $i=0,1,\ldots,l-j$ and then over
$\epsilon=\pm 1/2$ we get:
\begin{align*}
& \mathrm{ch}^F\!([e])
           =\sum_{l,j}\frac{(2j+1)(q^{2j+1}+q^{-2j-1})-[2][2j+1]}
           {[2l+2][2l+4][2j][2j+2]} \times \\ & 
           \qquad\times
           \sum_{2\epsilon=\pm 1}\Big\{(l-j+1)(q^{2l+2\epsilon(2j+1)+3}+q^{-2l-2\epsilon(2j+1)-3})
           -(q^{2\epsilon-1}+q^{-2\epsilon+1})[2]^{-1}[2(l-j+1)]\Big\} \\
           &\qquad=\sum_{l,j}\frac{(2j+1)(q^{2j+1}+q^{-2j-1})-[2][2j+1]}
           {[2l+2][2l+4][2j][2j+2]} \times \\ &
           \qquad\qquad\qquad\qquad\times
           \Big\{(l-j+1)(q^{2l+3}+q^{-2l-3})(q^{2j+1}+q^{-2j-1})
           -[2][2(l-j+1)]\Big\} \\
           &\qquad=:\sum_{l,j}f_{lj}(q)=:f(q) \;.
\end{align*}
We call $f_{lj}(q)$ the generic term of last series, explicitly written as
\begin{align*}
&f_{lj}(q)=(1-q^2)^4\,\frac{(2j+1)(1+q^{4j+2})-\frac{1+q^2}{1-q^2}(1-q^{4j+2})}
           {(1-q^{4l+4})(1-q^{4l+8})(1-q^{4j})(1-q^{4j+4})}
  \times \\ & \qquad\qquad\qquad\qquad\times
           q^{2l-1}\Big\{(l-j+1)(1+q^{4l+6})(1+q^{4j+2})
           -\tfrac{1+q^2}{1-q^2}\,q^2(q^{4j}-q^{4l+4})\Big\}\;,
\end{align*}
and consider it as a function of $q\in [\hspace{1pt}0,1[\,$. Notice that each
$f_{lj}(q)$ is a $C^\infty$ function of $q$ (they are rational functions whose
denominators never vanish for $0\leq q<1$). From the inequality
$$
0\leq f_{lj}(q)\leq 4(2j+1)q^{2l-1}
$$
we deduce (using the Weierstrass M-test) that the series is absolutely (hence uniformly)
convergent in each interval $[0,q_0]\subset [\hspace{1pt}0,1[\,$. Then, it converges to a
function $f(q)$ which is continuous in $[\hspace{1pt}0,1[\,$. 
Being the index of a Fredholm operator, $f(q)$ is integer valued in
$]\hspace{1pt}0,1[$; by continuity it is constant and can be computed in the limit
$q\to 0$. In this limit we have $f_{lj}(q)=2j(l-j+1)q^{2l-1}+O(q^{2l})$.
Thus, $f_{lj}(0)=\delta_{l,1/2}\delta_{j,1/2}$ and $\mathrm{ch}^F\!([e])=f(0)=1$.
\end{prova}

The next step is to define a spectral triple whose Fredholm module
is the one described in Proposition~\ref{prop:notr}.

\begin{prop}\label{prop:st}
Let $D$ be the (unbounded) operator on $\HH:=\HH_+\oplus\HH_-$ defined by
$$
D\ket{l,m_1,m_2;j}_\pm:=(l+\tfrac{3}{2})\ket{l,m_1,m_2;j}_\mp\;.
$$
Then, the datum $(\A(S^4_q),\HH,D,\gamma)$ is a $\Uq$-equivariant regular even
spectral triple of \mbox{metric} dimension $4$.
\end{prop}

\noindent{\bf Remark:} The operator $D$ is isospectral to the classical Dirac
operator on $S^4$ (whose spectrum has been computed in~\cite{Trau,CH}).
When $q=1$, this spectral triple becomes the canonical one associated to the spin
structure of $S^4$.

\begin{prova}
Clearly the representation of the algebra is even, $D$ is odd, with compact
resolvent and $4^+$-summable (being isospectral to the classical Dirac operator
on $S^4$).

Let $\delta$ be the unbounded derivation on $\B(\HH)$ defined by $\delta(T):=[|D|,T]$.
Each generator of $\A(S^4_q)$ is the sum of a finite number of weighted shifts;
each of these weighted shifts is a bounded operator (the coefficients are all
bounded by $1$) and is an eigenvector of $\delta$, i.e., if $T$ shifts the index
$l$ by $k$, then $\delta(T)=kT$. Thus, such weighted shifts are not
only bounded but also in the smooth domain of $\delta$, which we denote by
$\opz:=\bigcap_{j\in\N}\mathrm{dom}\,\delta^j$.
As a consequence $\A(S^4_q)\subset\opz$.

Recall that $[F,x_i]$ has coefficients decaying faster than $q^l$; thus
$|D|[F,x_i]$ is a matrix of rapid decay. In particular, $|D|[F,x_i]\in\op\subset\opz$.
The identity
\begin{equation}\label{eq:leib}
[D,x_i]=\delta(x_i)F+|D|[F,x_i]\;,
\end{equation}
tells us that $[D,x_i]$ is not only bounded but even in $\opz$ -- being the sum of two bounded operators contained in the $*$-algebra $\opz$. Then, $D$ defines a spectral triple and such a spectral triple is regular.

Finally, since $D$ is proportional to the identity in any irreducible subrepresentation
$V_l$ of $\Uq$, it commutes with all $h\in\Uq$ and it is equivariant.
\end{prova}

As a preparation for the study of the dimension spectrum in Sect.~\ref{sec:Sigma},
let us explicitly verify the $4$-summability of $D$. As one can easily check, the
dimension of $V_l$ is~\cite{CH}
$$
\dim V_l=\tfrac{2}{3}(l+\tfrac{5}{2})(l+\tfrac{3}{2})(l+\tfrac{1}{2}) \;.
$$
From this we get
$$
\tr(|D|^{-s})=\sum\nolimits_{l\in\N+\frac{1}{2}}2(l+\tfrac{3}{2})^{-s}\dim V_l
              =\tfrac{4}{3}\sum\nolimits_{n=1}^\infty (n^2-1)n^{-s+1} \;,
$$
where $n=l+\frac{3}{2}$ (and we added the term with $n=1$ since it is
identically zero). The above series is convergent in the right half-plane
$\{s\in\C\,|\,\mathrm{Re}\,s>4\}$, thus $D$ has metric dimension $4$.

Notice that $\tr(|D|^{-s})$ has meromorphic extension on $\C$ given by
\begin{equation}\label{eq:zeta}
\tr(|\D|^{-s})=\tfrac{4}{3}\big\{\zeta(s-3)-\zeta(s-1)\big\}\;,
\end{equation}
where $\zeta(s)$ is the Riemann zeta-function.
We recall that $\zeta(s)$ has a simple pole in $s=1$ as unique singularity
and that $\,\mathrm{Res}_{s=1}\zeta(s)=1$. 

\section{The dimension spectrum and residues}\label{sec:Sigma}
To compute the dimension spectrum we shall use a very simple representation of the algebra which differs -- in a sense which will be clear in Proposition~\ref{prop:app} -- from the chiral ones by a suitable ideal of operators. This is the class of operators,
\begin{equation}\label{eq:ideal}
\idJ:=\big\{T\in\opz\,\,\big|\,\, T |D|^{-p}\in\mathcal{L}^1(\HH)\,, \; \forall\;p>2\big\}\;.
\end{equation}

\begin{lemma}
The collection $\idJ$ is a two-sided ideal in $\opz$.
\end{lemma}
\begin{prova}
Clearly
$\idJ$ is a vector space: if $T_{1},T_{2}\in\idJ$, that is
$T_{1}|D|^{-p}\in\mathcal{L}^1(\HH)$, $T_{2}|D|^{-p}\in\mathcal{L}^1(\HH)$
for all $p>2$, then $T_{1}|D|^{-p}+T_{2}|D|^{-p}\in\mathcal{L}^1(\HH)$ for all $p>2$,
which means $T_1+T_2\in\idJ$.

That $\idJ$ is a left ideal is straightforward. 
Since $\mathcal{L}^1(\HH)$ is a two-sided ideal in $\B(\HH)$,
if $T_1\in\opz$ and $T_2\in\idJ$, for all $p>2$ we have that 
$T_1\cdot T_2|D|^{-p}$ is the product of a  bounded operator, $T_1$, with a trace class one, $T_2|D|^{-p}$, thus it is of trace class, and $T_1T_2\in\idJ$. 

From~Appendix B of \cite{CM} for any $p>0$, we know that the bounded operator $|D|^{-p}$ maps $\HH$ to $\HH^p:=\mathrm{dom}\,|D|^p$, that $T\in\opz\subset\mathrm{op}^0$
is a bounded operator $\HH^p\to\HH^p$, and finally that $|D|^p$ is
bounded from $\HH^p$ to $\HH$.  Thus, for $T\in\opz$, the product 
$|D|^pT|D|^{-p}$ is a bounded operator on $\HH$. 
Now, if $T_1\in\opz$ and $T_2\in\idJ$, for all $p>2$ we can write $T_2T_1|D|^{-p}=T_2 |D|^{-p} \cdot |D|^p T_1 |D|^{-p}$
as the product of a bounded operator, $|D|^{p}T_1|D|^{-p}$, with a trace class one, 
$T_2|D|^{-p}$; thus $T_2T_1|D|^{-p}$ is of trace class  so $T_2T_1\in\idJ$ and $\idJ$ is also a right ideal.
\end{prova}

\noindent
Clearly, if $T$ is of trace class, so is $|D|^{-p}T$ for any positive $p$,
and $\mathcal{L}^1(\HH)\subset\idJ$. Since $\op\subset\mathcal{L}^1(\HH)$,
smoothing operators belong to $\idJ$ as well. On the other hand, $\idJ$
is strictly bigger than $\mathcal{L}^1(\HH)$; indeed, 
the operator $L_q\in\B(\HH)$, given by
$$
L_q\ket{l,m_1,m_2;j}_\pm:=q^{j+\frac{1}{2}}\ket{l,m_1,m_2;j}_\pm\;,
$$
is not of trace class but belongs to $\idJ$, by the following proposition.
\begin{prop}\label{prop:jei}
For any $s\in\C$ with $\mathrm{Re}\,s>2$ one has that
$$
\zeta_{L_q}(s):=\sum_{l,j,m_1,m_2}(l+\tfrac{3}{2})^{-s}q^{j+\frac{1}{2}}=\frac{4q}{(1-q)^2}\,
\Big(\zeta(s-1)-\frac{1+q}{1-q}\,\zeta(s)\Big)+\textup{holomorphic function}\;,
$$
where $\zeta(s)$ is the Riemann zeta-function.
In particular, this means that $L_q$ belongs to the ideal $\idJ$.
Furthermore, since the series $\zeta_{L_q}(0)$ is divergent,
$L_q$ is not of trace class.
\end{prop}
\begin{prova}
Calling $n:=l+\frac{3}{2}$, $k:=j+\frac{1}{2}$, we have
$$
\zeta_{L_q}(s)=4\sum_{n=2}^\infty n^{-s}\sum_{k=1}^{n-1}k(n-k)q^k\;,
$$
We can sum starting
from $n=1$ and for $k=0,\ldots,n$ (we simply add zero terms) to get
$$
\zeta_{L_q}(s)=4\sum_{n=1}^\infty n^{-s}\big\{nq\partial_q-(q\partial_q)^2\big\}\sum_{k=0}^n q^k
    =4\sum_{n=1}^\infty n^{-s}\big\{nq\partial_q-(q\partial_q)^2\big\}\frac{1-q^{n+1}}{1-q}\;.
$$
Terms decaying as $q^n$ give a holomorphic function of $s$, thus modulo
holomorphic functions,
$$
\zeta_{L_q}(s)\sim 4\sum_{n=1}^\infty n^{-s}\left\{n\,\tfrac{q}{(1-q)^2}
-\tfrac{q(1+q)}{(1-q)^3}\right\} \;.
$$
The last series is summable for all $s$ with $\mathrm{Re}\,s>2$, and its sum
can be written in terms of the Riemann zeta-function as in the statement
of the proposition.
\end{prova}

\subsection{An approximated representation}\label{sec:HH}
Let $\hat{\HH}$ be a Hilbert space with orthonormal basis $\kkett{l,m_1,m_2;j}_\pm$
labelled by,
$$
l\in\tfrac{1}{2}\Z\;,\qquad
l+j\in\Z\;,\qquad
j+m_1\in\N\;,\qquad
l+\tfrac{1}{2}-j+m_2\in\N\;.
$$
Let $I$ be the labelling set of the Hilbert space $\HH_\pm$ as in Theorem
\ref{thm:cr}, and given by
$$
I:=\big\{(j,m_1,m_2,j)\;\big|\;l\in\N+\tfrac{1}{2}\;,\;\;j=\tfrac{1}{2},\tfrac{3}{2},...,l\;,
\;\;j-|m_1|\in\N\;,\;\;l+\tfrac{1}{2}-j-|m_2|\in\N\big\}\;.
$$
Notice that $I$ is the subset of labels of $\hat{\HH}$ satisfying
$l\in\N+\tfrac{1}{2}$, $m_1\leq j\leq l$ and $m_2\leq l+\tfrac{1}{2}-j$.
Define the inclusion $Q:\HH\to\hat{\HH}$ and the adjoint projection 
$P:\hat{\HH}\to\HH$ by,
\begin{align*}
Q\ket{l,m_1,m_2;j}_\pm &:=\kkett{l,m_1,m_2;j}_\pm 
\;\quad\textrm{for all}\;(l,m_1,m_2,j)\in I\;, \\
P\kkett{l,m_1,m_2;j}_\pm &:=\bigg\{\begin{array}{ll}
\ket{l,m_1,m_2;j}_\pm &\mathrm{if}\;(l,m_1,m_2,j)\in I\;, \\
0 &\mathrm{otherwise}\;.
\end{array}
\end{align*}
Clearly, $PQ=id_{\HH}$.
The Hilbert space $\hat{\HH}$ carries a bounded $*$-representation of the algebra $\A(SU_q(2))\otimes\A(S^2_q)$
defined by,
\begin{align*}
\alpha\kkett{l,m_1,m_2;j}_\pm &=\sqrt{1-q^{2(j+m_1+1)}}\kkett{l+\oh,m_1+\oh,m_2;j+\oh}_\pm \;,\\
\beta\kkett{l,m_1,m_2;j}_\pm  &=q^{j+m_1}\kkett{l+\oh,m_1-\oh,m_2;j+\oh}_\pm \;,\\
A\kkett{l,m_1,m_2;j}_\pm &=q^{l-j+m_2-\epsilon}\kkett{l,m_1,m_2;j}_\pm \;,\\
B\kkett{l,m_1,m_2;j}_\pm &=\sqrt{1-q^{2(l-j+m_2+2-\epsilon)}}\kkett{l+1,m_1,m_2+1;j}_\pm \;,
\end{align*}
where, as before, $\epsilon:=\frac{1}{2}(-1)^{l+\frac{1}{2}-j-m_2}$.
Composition of such a representation with the algebra embedding
$\A(S^4_q)\hookrightarrow\A(SU_q(2))\otimes\A(S^2_q)$ given in equation
(\ref{eq:embed}) results in a $*$-representation $\pi:\A(S^4_q)\to\B(\hat{\HH})$.
The sandwich $\tilde{\pi}(a):=P\pi(a)Q$ defines a $*$-linear map
$\tilde{\pi}:\A(S^4_q)\to\B(\HH)$.

\begin{prop}\label{prop:app}
With $\idJ$ the class of operators defined in Eq.~(\ref{eq:ideal}), one has that
the difference $a-\tilde{\pi}(a)\in\idJ$ for all $a\in\A(S^4_q)$.
\end{prop}
\begin{prova}
Define $\hat{\idJ}$ as the collection of bounded operators $T:\hat{\HH}\to\HH$ such that
$|D|^{-p}T$ is trace class for all $p>2$. Since trace class operators are a two sided ideal
in bounded operators, the space $\hat{\idJ}$ is stable when multiplied from the right
by bounded operators: $T_1\in\hat{\idJ}$ and $T_2\in\B(\hat{\HH})\;\Rightarrow\;T_1T_2\in\hat{\idJ}$.

Next, suppose that $a,b$ satisfy $a-\tilde{\pi}(a)\in\idJ$ and $b-\tilde{\pi}(b)\in\idJ$
and consider the following algebraic identity:
$$
ab-\tilde{\pi}(ab)=a\bigl\{b-\tilde{\pi}(b)\bigr\}+\bigl\{aP-P\pi(a)\bigr\}\pi(b)Q\;.
$$
Since $\idJ$ is a two-sided ideal in $\opz$, the first summand
is in $\idJ$. The stability of $\hat{\idJ}$ discussed above implies
that $\bigl\{aP-P\pi(a)\bigr\} \pi(b)\in\hat{\idJ}$, but if $T\in\hat{\idJ}$
clearly $TQ\in\idJ$. Hence the second summand in $\idJ$ too.
Thus, $ab-\tilde{\pi}(ab)\in\idJ$ whenever this property
holds for each of the operators $a,b$.
We conclude that it is enough to show that $a-\tilde{\pi}(a)\in\idJ$
when $a$ is a generator of $\A(S^4_q)$. By Proposition \ref{prop:jei},
this amounts to prove that the matrix elements of $a-\tilde{\pi}(a)$ 
are bounded in modulus by $q^{j+\frac{1}{2}}$.

Let us have a close look at the coefficients of $a\in\{x_i,x_i^*\}$ in Theorem \ref{thm:cr}.
Firstly, $A^+_{j,m_1}$, $B^+_{j,m_1}$, $B^-_{j,m_1}$, $q^{-2j}A^0_{j,m_1}$
and $q^{-2j}B^0_{j,m_1}$ are uniformly bounded by a constant, as one can see by writing explicitly the $q$-analogues in their expressions, getting:
\begin{align*}
A^+_{j,m_1} &=q^{j+m_1}(1-q^{4j+4})^{-1}\,
\sqrt{(1-q^{2(j+m_1+1)})(1-q^{2(j-m_1+1)})} \;, \\
q^{-2j}A^0_{j,m_1} &=(1-q^2)(1-q^{4j})^{-1}(1-q^{4j+4})^{-1} \,([2]q^{2(j+m_1)}-q^{4j+1}-q^{-1}) \;, \\
B^+_{j,m_1} &=(1-q^{4j+4})^{-1}\, \sqrt{(1-q^{2(j+m_1+1)})(1-q^{2(j+m_1+2)})}\;, \\
q^{-2j}B^0_{j,m_1} &=(1+q^2)q^{j+m_1+1}(1-q^{4j})^{-1}(1-q^{4j+4})^{-1} \, \sqrt{(1-q^{2(j-m_1)})(1-q^{2(j+m_1+1)})} \;, \\
B^-_{j,m_1} &=-q^{2(j+m_1)+1}(1-q^{4j})^{-1}\, 
\sqrt{(1-q^{2(j-m_1)})(1-q^{2(j-m_1-1)})} \;.
\end{align*}
Analogously, the coefficients
$q^{2j}H^0_{l,j,m_2}$, $C_{l,j,m_2}^0$ and $D_{l,j,m_2}^0$ are seen to be bounded by $q^l$ . Thus, modulo rapid decay matrices (i.e.~smoothing
operators),
\begin{subequations}\label{eq:apprep}
\begin{align}
x_0\ket{l,m_1,m_2;j} &\simeq A^+_{j,m_1}C^+_{l,j,m_2}\ket{l+1,m_1,m_2;j+1} \notag\\
                     &+A^+_{j,m_1}C^-_{l,j,m_2}\ket{l-1,m_1,m_2;j+1} \notag\\
                     &+A^0_{j,m_1}H^+_{l,j,m_2}\ket{l+1,m_1,m_2;j} \notag\\
                     &+A^0_{j,m_1}H^+_{l-1,j,m_2}\ket{l-1,m_1,m_2;j} \notag\\
                     &+A^+_{j-1,m_1}C^-_{l+1,j-1,m_2}\ket{l+1,m_1,m_2;j-1} \notag\\
                     &+A^+_{j-1,m_1}C^+_{l-1,j-1,m_2}\ket{l-1,m_1,m_2;j-1} \;, \\
x_1\ket{l,m_1,m_2;j} &\simeq B^+_{j,m_1}C^+_{l,j,m_2}\ket{l+1,m_1+1,m_2;j+1} \notag\\
                     &+B^+_{j,m_1}C^-_{l,j,m_2}\ket{l-1,m_1+1,m_2;j+1} \notag\\
                     &+B^0_{j,m_1}H^+_{l,j,m_2}\ket{l+1,m_1+1,m_2;j} \notag\\
                     &+B^0_{j,m_1}H^+_{l-1,j,m_2}\ket{l-1,m_1+1,m_2;j} \notag\\
                     &+B^-_{j,m_1}C^-_{l+1,j-1,m_2}\ket{l+1,m_1+1,m_2;j-1} \notag\\
                     &+B^-_{j,m_1}C^+_{l-1,j-1,m_2}\ket{l-1,m_1+1,m_2;j-1} \;, \\
x_2\ket{l,m_1,m_2;j} &\simeq D_{l,j,m_2}^+\ket{l+1,m_1,m_2+1;j} \notag\\
                     &+D_{l,j,m_2}^-\ket{l-1,m_1,m_2+1;j} \;.
\end{align}
\end{subequations}
Since modulo smoothing operators the representations are the same we are omitting the label `$\pm$' in the vector basis.
Furthermore, using the inequalities
\begin{equation}\label{eq:ineqJ}
0\leq (1-qu)^{-1}-1\leq q (1-q)^{-1} \,u\;,\qquad
0\leq 1-(1-u)^{\frac{1}{2}}\leq u\;,
\end{equation}
which are valid when $0\leq u\leq 1$,
we prove that modulo terms bounded by $q^l$, one has
\begin{subequations}\label{eq:coefsim}
\begin{align}
C^+_{l,j,m_2} &\simeq -q^{l-j+m_2-\epsilon}\sqrt{1-q^{2(l+j+m_2+3+\epsilon)}} \;, \\
C^-_{l,j,m_2} &\simeq -q^{l+j+m_2+1+\epsilon}\sqrt{1-q^{2(l-j+m_2-\epsilon)}} \;, \\
H^+_{l,j,m_2} &\simeq q^{l+m_2+1}\sqrt{q^{2\epsilon(2j+1)}-q^{2(l+m_2+2)}} \;, \\
D_{l,j,m_2}^+ &\simeq \sqrt{1-q^{2(l+j+m_2+3+\epsilon)}}\sqrt{1-q^{2(l-j+m_2+2-\epsilon)}} \;, \\
D_{l,j,m_2}^- &\simeq -q^{2(l+m_2)+3} \;.
\end{align}
\end{subequations}
Up to now, we neglected only smoothing contributions (the above approximation
will be needed when dealing with the real structure later on). We use again
(\ref{eq:ineqJ}) to get a rougher approximation by neglecting terms bounded by $q^j$. This yields
\begin{subequations}\label{eq:get}
\begin{align}
A^+_{j,m_1} &\simeq \tilde{A}^+_{j,m_1}:=q^{j+m_1}\sqrt{1-q^{2(j+m_1+1)}} \;, \\
A^0_{j,m_1}H^+_{l,j,m_2}  &\simeq 0 \;, \\
B^+_{j,m_1} &\simeq \tilde{B}^+_{j,m_1}:=\sqrt{(1-q^{2\smash[t]{(j+m_1+1)}})(1-q^{2\smash[t]{(j+m_1+2)}})} \;, \\
B^0_{j,m_1}H^+_{l,j,m_2}  &\simeq 0 \;, \\
B^-_{j,m_1} &\simeq \tilde{B}^-_{j,m_1}:=-q^{2(j+m_1)+1} \;,\\
C^+_{l,j,m_2} &\simeq \tilde{C}_{l,j,m_2}:=-q^{l-j+m_2-\epsilon} \;, \\
C^-_{l,j,m_2} &\simeq 0 \;, \\
D_{l,j,m_2}^+ &\simeq \tilde{D}_{l,j,m_2}:=\sqrt{1-q^{2(l-j+m_2+2-\epsilon)}} \;, \\
D_{l,j,m_2}^- &\simeq 0 \;.
\end{align}
\end{subequations}
Plugging these coefficients in the equations for the $x_i$'s we see that,
modulo operators in the ideal $\idJ$, we get
\begin{align*}
x_0\ket{l,m_1,m_2;j} &\simeq \tilde{A}^+_{j,m_1}\tilde{C}_{l,j,m_2}\ket{l+1,m_1,m_2;j+1} \\
                     &\quad +\tilde{A}^+_{j-1,m_1}\tilde{C}_{l-1,j-1,m_2}\ket{l-1,m_1,m_2;j-1} \\
                     &=-P(\alpha\beta A)Q\ket{l,m_1,m_2;j} - P(\beta^*\alpha^* A)Q\ket{l,m_1,m_2;j} \;, \\
x_1\ket{l,m_1,m_2;j} &\simeq \tilde{B}^+_{j,m_1}\tilde{C}_{l,j,m_2}\ket{l+1,m_1+1,m_2;j+1} \\
                     &\quad +\tilde{B}^-_{j,m_1}\tilde{C}_{l-1,j-1,m_2}\ket{l-1,m_1+1,m_2;j-1} \\
                     &=-P(\alpha^2 A)Q\ket{l,m_1,m_2;j} + P(q(\beta^*)^2 A)Q\ket{l,m_1,m_2;j} \;, \\
x_2\ket{l,m_1,m_2;j} &\simeq \tilde{D}_{l,j,m_2}\ket{l+1,m_1,m_2+1;j} \\
                     &=PB\,Q \ket{l+1,m_1,m_2+1;j} \;.
\end{align*}
The observation that
$$
-P(\alpha\beta+\beta^*\alpha^*)AQ=\tilde{\pi}(x_0)\;,\qquad
P(-\alpha^2+q(\beta^*)^2)AQ=\tilde{\pi}(x_1)\;,\qquad
PBQ=\tilde{\pi}(x_2)\;,
$$
concludes the proof.
\end{prova}

\subsection{The dimension spectrum and the top residue}
The approximation modulo $\idJ$ allows considerable simplifications when getting information on the part of
the dimension spectrum contained in the half plane $\mathrm{Re}\,s>2$.
To study the part of the dimension spectrum in the left half plane
$\mathrm{Re}\,s\leq 2$ would require a less drastic approximation which we are lacking at the moment.

\begin{prop}
In the region $\mathrm{Re}\,s>2$ the dimension spectrum $\Sigma$ of the spectral triple $(\A(S^4_q),\HH,D,\gamma)$ given in Proposition~\ref{prop:st} consists of the
two points $\{3,4\}$, which are simple poles of the zeta-functions. The top residue
coincides with the integral on the subspace of classical points of $S^4_q$, that is
\begin{equation}\label{eq:topres}
\nint a|D|^{-4}=\frac{2}{3\pi}\int_0^{2\pi}\!\!\sigma(a)(\theta)\de\theta\;,
\end{equation}
with $\sigma:\A(S^4_q)\to\A(S^1)$ the $*$-algebra morphism defined by
$\sigma(x_0)=\sigma(x_1)=0$ and $\sigma(x_2)=u$, where $u$, given by
$u(\theta):=e^{i\theta}$, is the unitary generator of $\A(S^1)$.
\end{prop}
\begin{prova}
Let $\Psi^0$ be the $*$-algebra generated by $\A(S^4_q)$, by $[D,\A(S^4_q)]$
and by iterated applications of the derivation $\delta$. Let $\mathfrak{A}\subset
\A(SU_q(2))\otimes\A(S^2_q)\otimes\mathrm{Mat}_2(\C)$ be the $*$-algebra generated
by $\alpha$, $\beta$, $\alpha^*$, $\beta^*$, $A$, $B$, $B^*$ and $F$.
By Proposition~\ref{prop:app} there is an inclusion $\A(S^4_q)\subset P\mathfrak{A}Q+\idJ$.

A linear basis for $\mathfrak{A}$ is given by,
\begin{equation}\label{eq:linb}
T:=\alpha^{k_1}\beta^{n_1}(\beta^*)^{n_2}A^{n_3}B^{k_2}F^h\;,
\end{equation}
where $h\in\{0,1\}$, $n_i\in\N$, $k_i\in\Z$ and
with the notation $\alpha^{k_1}:=(\alpha^*)^{-k_1}$ if $k_1<0$ and
$B^{k_2}:=(B^*)^{-k_2}$ if $k_2<0$. For this operator,
$$
\delta(PTQ)=\big(\tfrac{1}{2}(k_1+n_1-n_2)+k_2\big)PTQ\qquad\mathrm{and}\qquad
[D,PTQ]=\delta(PTQ)F\;.
$$
Thus, $P\mathfrak{A}Q$ is invariant under application of $\delta$ and $[D,\,.\,]$
and hence $\Psi^0\subset P\mathfrak{A}Q+\idJ$.

For the part of the dimension spectrum in the right half
plane $\mathrm{Re}\,s>2$, we can neglect $\idJ$ and consider only the singularities
of zeta-functions associated with elements in $P\mathfrak{A}Q$. By linearity of
the zeta-functions, it is enough to consider the generic basis element in
Eq.~(\ref{eq:linb}).

Such a $T$ shifts $l$ by $\frac{1}{2}(k_1+n_1-n_2)+k_2$, $m_1$ by $\frac{1}{2}(k_1-n_1+n_2)$,
$m_2$ by $k_2$, $j$ by $\frac{1}{2}(k_1+n_1-n_2)$, and flips the chirality if $h=1$.
Thus it is off-diagonal unless $h=k_i=0$ and $n_1=n_2$. The zeta-function associated with
a bounded off-diagonal operator is identically zero in the half-plane $\mathrm{Re}\,z>4$,
and so is its holomorphic extension to the entire complex plane. It remains to consider the cases $T=P(\beta\beta^*)^kA^nQ$, with $n,k\in\N$.

If $n$ and $k$ are both different from zero, one finds
$$
\zeta_T(s)=2\sum_{l,j,m_1,m_2} (l+\tfrac{3}{2})^{-s}q^{n(l-j+m_2-\epsilon)+2k(j+m_1)}
 =2\sum_{l,j,m_2} (l+\tfrac{3}{2})^{-s}q^{n(l-j+m_2-\epsilon)}\frac{1-q^{2k(2j+1)}}{1-q^{2k}}\;.
$$
For $\epsilon$ fixed, set $2i:=l-\epsilon-j+m_2\in\{0,2,\ldots,2(l-j)\}$. Then,
\begin{align*}
\zeta_T(s)&=2\sum_{l,j} (l+\tfrac{3}{2})^{-s}\frac{1-q^{2k(2j+1)}}{1-q^{2k}}
\sum_{\epsilon=\pm 1/2}\sum_{i=0}^{l-j}q^{2ni}=
4\sum_{l,j} (l+\tfrac{3}{2})^{-s}\frac{1-q^{2k(2j+1)}}{1-q^{2k}}
\frac{1-q^{2n(l-j+1)}}{1-q^{2n}} \\
&=4\zeta(s-1)-4\,\frac{1+(1-q^{4k})^{-1}+(1-q^{2n})^{-1}}{(1-q^{2k})(1-q^{2n})}\,\zeta(s)
  +\textrm{holomorphic function}\;, 
\end{align*}
which has meromorphic extension on $\C$ with simple pole in $s=\{1,2\}$.

If $n=0$ and $k\neq 0$,
\begin{align*}
\zeta_T(s) &=4\sum_{l,j}(l+\tfrac{3}{2})^{-s}(l-j+1)\frac{1-q^{2k(2j+1)}}{1-q^{2k}} \\
  &=\tfrac{4}{1-q^{2k}}\left(\tfrac{1}{2}\,\zeta(s-2)-\big(\tfrac{1}{2}+\tfrac{1}{1-q^{4k}}\big)
  \zeta(s-1)+\tfrac{q^{4k}}{(1-q^{4k})^2\log q^{4k}}\,\zeta(s)\right)
  +\textrm{hol.~function}\;,
\end{align*}
which has meromorphic extention on $\C$ with simple pole in $s=\{1,2,3\}$.

If $n\neq 0$ and $k=0$, 
\begin{align*}
\zeta_T(s)&=4\sum_{l,j} (l+\tfrac{3}{2})^{-s}(2j+1)\frac{1-q^{2n(l-j+1)}}{1-q^{2n}} \\
  &=\tfrac{4}{1-q^{2n}}\left\{\zeta(s-2)-\Big(1+\tfrac{2q^{2n}}{1-q^{2n}}\Big)\zeta(s-1)
  +\tfrac{2q^{2n}}{1-q^{2n}}\Big(1+\tfrac{q^{2n}}{(1-q^{2n})\log q^{2n}}\Big)\zeta(s)\right\}
  +\textrm{hol.~fun.}\;, 
\end{align*}
which has meromorphic extention on $\C$ with simple pole in $s=\{1,2,3\}$.

Finally, if both $n$ and $k$ are zero we get (cf.~Eq.~(\ref{eq:zeta})),
$$
\zeta_T(s)=\tfrac{4}{3}\big\{\zeta(s-3)-\zeta(s-1)\big\}\;,
$$
and this is meromorphic with simple poles in $\{2,4\}$. Thus, the part of the
dimension spectrum in the region $\mathrm{Re}\,s>2$ consists at
most of the two points $\{3,4\}$ and both are simple poles.

Since we have considered the enlarged algebra $P\mathfrak{A}Q+\idJ$, it suffices
to prove that there exists an $a\in\Psi^0$ whose zeta-function is
singular in both points $s=3$ and $s=4$. We take $a=x_2x_2^*$. From the definition
$$
\tilde{\pi}(x_2x_2^*)\ket{l,m_1,m_2;j}_\pm=
(1-q^{2(l-\epsilon-j+m_2)})\ket{l,m_1,m_2;j}_\pm\;.
$$
Then, modulo functions that are holomorphic when $\mathrm{Re}\,s>2$, we have
$$
\zeta_{x_2x_2^*}(s)\sim\zeta_{\tilde{\pi}(x_2x_2^*)}(s)= \zeta_1(s)
-2\sum_{l,j,m_1,m_2}(l+\tfrac{3}{2})^{-s}q^{2(l-\epsilon-j+m_2)}
\sim\tfrac{4}{3}\,\zeta(s-3)-\tfrac{4}{1-q^4}\,\zeta(s-2)\;.
$$
This proves the first part of the proposition,
that is $\Sigma\cap\{\mathrm{Re}\,s>2\}=\{3,4\}$.

The proof of Eq.~(\ref{eq:topres}) is based on the observation that the residue
in $s=4$ of $\zeta_T$, for $T$ a basis element of $P\mathfrak{A}Q$, is zero unless
$T=1$. That is, it depends only on the image of $T$ under the map sending $\beta,A$ and $F$ to $0$ while $\alpha\mapsto e^{i\phi}$ and
$B\mapsto e^{i\theta}$. Composing this map with $\tilde{\pi}$ we get the morphism $\sigma:\A(S^4_q)\to \A(S^1)$ of the proposition and that
$$
\nint a|D|^{-4}\propto\int_0^{2\pi}\!\!\sigma(a)\de\theta\;.
$$
The equality $\,\int\mkern-16mu- \,|D|^{-4}=\frac{4}{3}\,$ fixes the proportionality
constant.
\end{prova}

\section{Reality and first order conditions}\label{sec:real}
Classically, if $(\A(M),\HH,D,\gamma)$ is the canonical spectral triple associated with a
$4$-dimensional spin manifold $M$, there exists an antilinear isometry $J$ on $\HH$,
named the \emph{real structure}, satisfying the following compatibility condition
\begin{equation}\label{eq:noA}
J^2=-1\;,\quad\qquad
J\gamma=\gamma J\;,\quad\qquad
JD=DJ\;.
\end{equation}
There are also two additional conditions involving the coordinate algebra $\A(M)$:
\begin{equation}\label{eq:further}
[a,JbJ^{-1}]=0\;,\qquad[[D,a],JbJ^{-1}]=0\;,\qquad\forall\;a,b\in\A(M)\;.
\end{equation}
The real structure on $S^4$ is equivariant and equivariance is sufficient to
determine $J$. 

In the deformed situation one has to be careful on how to implement equivariance.
Let us start with the working hypothesis that equivariance for $J$ is the
requirement that it satisfies $Jh=S(h)^*J$ for all $h\in\Uq$. 
Then, consider the Casimir operator $\mathcal{C}_1$ given in equation
(\ref{eq:Cdef}). This operator commutes with $J$ since $S(\mathcal{C}_1)^*=\mathcal{C}_1$ and from its expression, $\mathcal{C}_1\ket{l,m_1,m_2;j}=(q^{2j+1}+q^{-2j-1})\ket{l,m_1,m_2;j}$,
we conclude that $J$ leaves the index $j$ invariant.
Compatibility with $\gamma$ and $D$ in Eq.~(\ref{eq:noA}) and equivariance
with respect to $h=K_1$ and $h'=K_2$ yields
$$
J\ket{l,m_1,m_2;j}_\pm=c_\pm(l,m_1,m_2;j)\ket{l,-m_1,-m_2;j}_\pm\;,
$$
with some constants  $c_\pm$  to be determined. Equivariance with respect to $h=E_1$
implies
$$
c_\pm(l,m_1,m_2;j)=(-1)^{m_1+1/2}q^{m_1}c_\pm(l,m_2;j)\;.
$$
For $h=E_2$, looking at the piece diagonal in $j$ we deduce that the dependence on $m_2$
is through a factor $q^{3m_2}$; and looking at the piece shifting $j$ by $\pm 1$ we conclude that
$$
c_\pm(l,m_1,m_2;j)=(-1)^{j+m_1}q^{m_1+3m_2}c_\pm(l)\;.
$$
Such an operator $J$  cannot be antiunitary unless $q=1$. At $q=1$ the antiunitarity
condition requires that $c_\pm(l)\in U(1)$ and modulo a unitary equivalence we can choose
$c_\pm(l)=i^{2l+1}$. In conclusion for $q=1$ the operator
\begin{equation}\label{eq:J}
J\ket{l,m_1,m_2;j}_\pm=i^{2l+1}(-1)^{j+m_1}\ket{l,-m_1,-m_2;j}_\pm\;,
\end{equation}
is \emph{the} real structure on $S^4$ (modulo a unitary equivalence). 

For $q\neq 1$ we
keep (\ref{eq:J}) as the real structure and notice that conditions \eqref{eq:noA} are
satisfied, but $J$ no longer satisfies the requirement $Jh=S(h)^*J$ for all $h\in\Uq$. 
Nevertheless, $J$ is the
antiunitary part of an  antilinear operator $T$ that has this property. The antilinear operator $T$ defined by 
$$
T\ket{l,m_1,m_2;j}_\pm=i^{2l+1}(-1)^{j+m_1}q^{m_1+3m_2}\ket{l,-m_1,-m_2;j}_\pm\;,
$$
has $J$ in \eqref{eq:J} as the antiunitary part and it is \emph{equivariant}, i.e. it is such that $Th=S(h)^*T$ for all $h\in\Uq$.

Next, we turn to the conditions (\ref{eq:further}). 
In parallel with the cases of the manifold of $SU_q(2)$ in \cite{DLSvSV} and
of Podle{\'s} spheres in \cite{DLPS,DDLW}, once again we need to modify
them.
For instance, the commutator $[x_2,Jx_2J]$ is not zero, as one can see
by computing the matrix element
\begin{align*}
f(l,j,m_2)&:=\pm \inner{l+1,m_1,m_2;j\big|[x_2,Jx_2J]\big|l,m_1,m_2;j}_\pm \\
&=D_{l+1,j,m_2-1}^0D_{l,j,-m_2}^+-D_{l,j,m_2-1}^+D_{l,j,-m_2}^0
+D_{l+1,j,-m_2-1}^0D_{l,j,m_2}^+-D_{l,j,-m_2-1}^+D_{l,j,m_2}^0 \;,
\end{align*}
which for $j=\frac{1}{2}$ and $m_2=l$ is
$$
f(l,\tfrac{1}{2},l)=-q^{-l-4}(1-q^2)^2[2]\frac{(q^{l-1}+q^{-l+1})
\sqrt{[2l+3]}\,[l+1][l+2][l+3]}{[2l+2][2l+4]^2[2l+6]}\neq 0\;.
$$
It is relatively easy to prove that the two conditions are satisfied
modulo the ideal $\idJ$. It is much more cumbersome computationally to
show that they are in fact satisfied modulo the smaller ideal of
smoothing operators.

\begin{prop}\label{prop:real}
Let $J$ be the antilinear isometry given by (\ref{eq:J}). Then,
$$
[a,JbJ]\in\idJ\;,\quad[[D,a],JbJ]\in\idJ\;,\qquad\forall\;a,b\in\A(S^4_q)\;.
$$
\end{prop}
\begin{prova}
We lift $J$ and $D$ to the Hilbert space $\hat{\HH}$
defined in Sect.~\ref{sec:HH}, as follows:
\begin{align*}
\hat{J}\!\kkett{l,m_1,m_2;j}_\pm &=i^{2l+1}(-1)^{j+m_1}\!\kkett{l,-m_1,-m_2;j}_\pm \;, \\
\hat{D}\!\kkett{l,m_1,m_2;j}_\pm &=(l+\tfrac{3}{2})\!\kkett{l,m_1,m_2;j}_\mp \;.
\end{align*}
Notice that $\hat{J}^2=-1$ on $\hat{\HH}$ (thanks to the phase $i^{2l+1}$ that is
irrelevant when restricted to $\HH$).

Let now $\{\alpha,\beta,\alpha^*,\beta^*,A,B,B^*\}$ be the operators defined in Section
\ref{sec:HH}, generators of the algebra $\A(SU_q(2))\otimes\A(S^2_q)$.
Due to Proposition~\ref{prop:app} it is enough to prove that for all
pairs $(a,b)$ of such generators,
the commutators $[a,\hat{J}b\hat{J}]$ and $[[\hat{D},a],\hat{J}b\hat{J}]$
are weighted shifts with weight which are bounded by $q^{2j}$. From 
$$
[\hat{D},\alpha]=\tfrac{1}{2}\,\alpha\hat{F}\;,\quad
[\hat{D},\beta]=\tfrac{1}{2}\,\beta\hat{F}\;,\quad
[\hat{D},A]=0\;,\quad
[\hat{D},B]=B\hat{F}\;,
$$
the condition on $[[\hat{D},a],\hat{J}b\hat{J}]$ follows
from the same condition on $[a,\hat{J}b\hat{J}]$,
and we have to compute only the latter commutators.

Since $[a,\hat{J}b^*\hat{J}]=-[a^*,\hat{J}b\hat{J}]^*$ and 
$[b,\hat{J}a\hat{J}]=\hat{J}[a,\hat{J}b\hat{J}]\hat{J}$,
we have to check the $16$ combinations in the following table.

\begin{center}\begin{tabular}{c|ccccccc}
\hline\hline $b \backslash a$ &
           $\alpha$ & $\alpha^*$ & $\beta$ & $\beta^*$ & $A$ & $B$ & $B^*$ \\
\hline
$\alpha$   & $\bullet$ & $\bullet$ & $\times$ & $\times$ & $\bullet$ & $\bullet$ & $\bullet$ \\
$\beta$    &        &          & $\bullet$ & $\times$ & $\bullet$ & $\bullet$ & $\bullet$ \\
$B$        &        &          &       &         & $\bullet$ & $\bullet$ & $\bullet$ \\
$A$        &        &          &       &         & $\bullet$ &   & \\
\hline\hline
\end{tabular}
\end{center}

By direct computations one shows that bullets in the table correspond to
vanishing commutators. On the other hand, the commutators corresponding to
the crosses in the table are given, on the subspace with $j-|m_1|\in\N$, by
\begin{align*}
[\beta^*,\hat{J}\alpha\hat{J}]\kkett{l,m_1,m_2;j}_\pm &=q^{j+m_1}\Big\{\sqrt{1-q^{2(j-m_1+1)}}-
\sqrt{1-q^{2(j-m_1)}}\Big\}\kkett{l,m_1,m_2;j}_\pm \\
[\beta,\hat{J}\alpha\hat{J}]\kkett{l,m_1,m_2;j}_\pm &=-[\beta^*,\hat{J}\alpha\hat{J}]\kkett{l+1,m_1-1,m_2;j+1}_\pm\;,\\
[\beta^*,\hat{J}\beta\hat{J}]\kkett{l,m_1,m_2;j}_\pm &=-[2]q^{2j}\kkett{l,m_1+1,m_2;j}_\pm \;.
\end{align*}
Since $1-u\leq\sqrt{1-u}\leq 1$ for all $u\in[0,1]$, we have that
$$
0\leq q^{j+m_1}\Big\{\sqrt{1-q^{2(j-m_1+1)}}-\sqrt{1-q^{2(j-m_1)}}\Big\}
\leq q^{j+m_1}(1-1+q^{2(j-m_1)})\leq q^{2j}\;.
$$
Then, all three non-zero commutators are weighted shifts with weights bounded by
$q^{2j}$. 
\end{prova}

\begin{prop}
Let $J$ be the antilinear isometry given by (\ref{eq:J}). Then,
$$
[a,JbJ]\in\op\;,\quad[[D,a],JbJ]\in\op\;,\qquad\forall\;a,b\in\A(S^4_q)\;.
$$
\end{prop}

\begin{prova}
By Leibniz rule, it is sufficient to prove the statement
when $a$ and $b$ are generators of the algebra. By (\ref{eq:leib}),
$[D,a]-\delta(a)F$ is a smoothing operator.
Thus,
it is enough to show that
\begin{equation}\label{eq:aJbJ}
[a,JbJ]\in\op\;,\qquad
[\delta(a),JbJ]\in\op\;,
\end{equation}
for any pair $(a,b)$ of generators.
From
$$
[b,JaJ]=J[a,JbJ]J \;,\qquad
[\delta(b),JaJ]=-J[\delta(a),JbJ]J+\delta([b,JaJ]) \;,
$$
it follows that if (\ref{eq:aJbJ}) is satisfied for a particular
pair $(a,b)$, then it is satisfied for $(b,a)$ too. From
\begin{equation}\label{eq:symstar}
[a^*,Jb^*J]=-[a,JbJ]^*\;,\qquad
[\delta(a^*),b^*]=[\delta(a),b]^* \;,
\end{equation}
we see that if (\ref{eq:aJbJ}) is satisfied for a
pair $(a,b)$, then it is satisfied for $(a^*,b^*)$ too. 
With these symmetries we need to check only the following
$9$ cases out of $25$:

\begin{center}\begin{tabular}{c|ccccc}
\hline\hline $b \backslash a$ &
           $x_0$ & $x_1$ & $x_1^*$ & $x_2$ & $x_2^*$ \\
\hline
$x_0$   & $\bullet$ & $\bullet$ && $\bullet$ \\
$x_1$  && $\bullet$ & $\bullet$ & $\bullet$ & $\bullet$ \\
$x_2$&&&& $\bullet$ & $\bullet$ \\
\hline\hline
\end{tabular}
\end{center}

\noindent
From Eqs.~(\ref{eq:apprep}) and (\ref{eq:coefsim}) we see that modulo smoothing operators
\begin{align*}
\oh\bigl\{x_0+\delta(x_0)\bigr\}\ket{l,m_1,m_2;j} 
        &\simeq \boxed{ A^+_{j,m_1}\hat{C}^+_{l,j,m_2}\ket{l+1,m_1,m_2;j+1} } \\
        &+q^{-2j}A^0_{j,m_1}\hat{H}^+_{l,j,m_2}\ket{l+1,m_1,m_2;j} \\
        &+A^+_{j-1,m_1}\hat{C}^-_{l+1,j-1,m_2}\ket{l+1,m_1,m_2;j-1} \;,\\
\oh\bigl\{x_0-\delta(x_0)\bigr\}\ket{l,m_1,m_2;j}
        &\simeq A^+_{j,m_1}\hat{C}^-_{l,j,m_2}\ket{l-1,m_1,m_2;j+1} \\
        &+q^{-2j}A^0_{j,m_1}\hat{H}^+_{l-1,j,m_2}\ket{l-1,m_1,m_2;j} \\
        &+\boxed{ A^+_{j-1,m_1}\hat{C}^+_{l-1,j-1,m_2}\ket{l-1,m_1,m_2;j-1} } \;, \\
\oh\bigl\{x_1+\delta(x_1)\bigr\}\ket{l,m_1,m_2;j}
        &\simeq \boxed{ B^+_{j,m_1}\hat{C}^+_{l,j,m_2}\ket{l+1,m_1+1,m_2;j+1} } \\
        &+q^{-2j}B^0_{j,m_1}\hat{H}^+_{l,j,m_2}\ket{l+1,m_1+1,m_2;j} \\
        &+B^-_{j,m_1}\hat{C}^-_{l+1,j-1,m_2}\ket{l+1,m_1+1,m_2;j-1} \;,\\
\oh\bigl\{x_1-\delta(x_1)\bigr\}\ket{l,m_1,m_2;j}
        &\simeq B^+_{j,m_1}\hat{C}^-_{l,j,m_2}\ket{l-1,m_1+1,m_2;j+1} \\
        &+q^{-2j}B^0_{j,m_1}\hat{H}^+_{l-1,j,m_2}\ket{l-1,m_1+1,m_2;j} \\
        &+\boxed{ B^-_{j,m_1}\hat{C}^+_{l-1,j-1,m_2}\ket{l-1,m_1+1,m_2;j-1} } \;, \\
\oh\bigl\{x_2+\delta(x_2)\bigr\}\ket{l,m_1,m_2;j}
        &\simeq \boxed{\boxed{ \hat{D}_{l,j,m_2}^+\ket{l+1,m_1,m_2+1;j} }} \;,\\
\oh\bigl\{x_2-\delta(x_2)\bigr\}\ket{l,m_1,m_2;j}
        &\simeq\hat{D}_{l,j,m_2}^-\ket{l-1,m_1,m_2+1;j} \;.
\end{align*}
where
\begin{subequations}\label{eq:Dhat}
\begin{align}
\hat{C}^+_{l,j,m_2} &=-q^{l-j+m_2-\epsilon}\sqrt{1-q^{2(l+j+m_2+3+\epsilon)}} \;, \\
\hat{C}^-_{l,j,m_2} &=-q^{l+j+m_2+1+\epsilon}\sqrt{1-q^{2(l-j+m_2-\epsilon)}} \;, \\
\hat{H}^+_{l,j,m_2} &=q^{2j}q^{l+m_2+1}\sqrt{q^{2\epsilon(2j+1)}-q^{2(l+m_2+2)}} \;, \\
\hat{D}_{l,j,m_2}^+ &=\sqrt{1-q^{2(l+j+m_2+3+\epsilon)}}\sqrt{1-q^{2(l-j+m_2+2-\epsilon)}} \;,\\
\hat{D}_{l,j,m_2}^- &=-q^{2(l+m_2)+3} \;.
\end{align}
\end{subequations}
We have divided the terms in three classes, which need to be analysed
separately.

All terms $T$ which are not `boxed' have coefficients which are uniformly
bounded by $q^{l+m_2}$; since the conjugation with $J$ changes the sign of
the labels $m_1,m_2$, for such $T$'s, the coefficients of $JTJ$ are
uniformly bounded by $q^{l-m_2}$.
They give products (and so commutators) with coefficients bounded by
$q^{l+m_2}q^{l-m_2}=q^{2l}$, and so (these products) are smoothing operators.

Analogously, the coefficients of single-boxed terms are bounded by $q^{l-j+m_2}$,
and become smoothing when multiplied by the $J$-conjugated of non-boxed
terms (as $q^{l-j+m_2}q^{l-m_2}\leq q^l$), and viceversa for the
product of a non-boxed term with the $J$-conjugated of a single-boxed
one ($q^{l+m_2}q^{l-j-m_2}\leq q^l$).

Next we consider pairs of single-boxed terms.
A closer look at the single-boxed terms in $x_0\pm\delta(x_0)$ and $x_1-\delta(x_1)$
(and then $x_1^*+\delta(x_1^*)$) shows that they have coefficients bounded by $q^{l+m_1+m_2}$,
and become smoothing when multiplied by the $J$-conjugated of one of them
($q^{l+m_1+m_2}q^{l-m_1-m_2}=q^{2l}$).
Last single-boxed term is the one in $x_1+\delta(x_1)$ (and $x_1^*-\delta(x_1^*)$).
The relevant terms for the commutators involving them are
\begin{align*}
&\oh[x_1+\delta(x_1),Jx_0J]\ket{l,m_1,m_2;j} \\
&\qquad \simeq\bigl\{A^+_{j+1,-m_1-1}\hat{C}^+_{l+1,j+1,-m_2}B^+_{j,m_1}\hat{C}^+_{l,j,m_2} + \\ &\qquad\qquad\qquad\qquad
         -A^+_{j,-m_1}\hat{C}^+_{l,j,-m_2}B^+_{j+1,m_1}\hat{C}^+_{l+1,j+1,m_2}\bigr\}\ket{l+2,m_1+1,m_2;j+2} \\
&\qquad\qquad +\bigl\{A^+_{j,-m_1-1}\hat{C}^+_{l,j,-m_2}B^+_{j,m_1}\hat{C}^+_{l,j,m_2} + \\ &\qquad \qquad\qquad\qquad
         -A^+_{j-1,-m_1}\hat{C}^+_{l-1,j-1,-m_2}B^+_{j-1,m_1}\hat{C}^+_{l-1,j-1,m_2}\bigr\}\ket{l,m_1+1,m_2;j} \;, \\
&\oh[x_1+\delta(x_1),Jx_1J]\ket{l,m_1,m_2;j} \\
&\qquad \simeq\bigl\{B^+_{j+1,m_1-1}\hat{C}^+_{l+1,j+1,m_2}B^+_{j,-m_1}\hat{C}^+_{l,j,-m_2} + \\ &\qquad\qquad\qquad\qquad
         -B^+_{j+1,-m_1-1}\hat{C}^+_{l+1,j+1,-m_2}B^+_{j,m_1}\hat{C}^+_{l,j,m_2}\bigr\}\ket{l+2,m_1,m_2;j+2} \\
&\qquad\qquad +\bigl\{B^+_{j-1,m_1-1}\hat{C}^+_{l-1,j-1,m_2}B^-_{j,-m_1}\hat{C}^+_{l-1,j-1,-m_2} + \\ &\qquad\qquad\qquad\qquad
         -B^-_{j+1,-m_1-1}\hat{C}^+_{l,j,-m_2}B^+_{j,m_1}\hat{C}^+_{l,j,m_2}\bigr\}\ket{l,m_1,m_2;j} \;,\\
&\oh[x_1-\delta(x_1),Jx_1J]\ket{l,m_1,m_2;j} \\
&\qquad \simeq\bigl\{B^-_{j+1,m_1-1}\hat{C}^+_{l,j,m_2}B^+_{j,-m_1}\hat{C}^+_{l,j,-m_2} + \\ &\qquad\qquad\qquad\qquad
         -B^-_{j,m_1}\hat{C}^+_{l-1,j-1,m_2}B^+_{j-1,-m_1-1}\hat{C}^+_{l-1,j-1,-m_2}\bigr\}\ket{l,m_1,m_2;j} \;,\\
&\oh[x_1^*+\delta(x_1^*),Jx_1J]\ket{l,m_1,m_2;j} \\
&\qquad\simeq \bigl\{B^-_{j+2,m_1-2}\hat{C}^+_{l+1,j+1,m_2}B^+_{j,-m_1}\hat{C}^+_{l,j,-m_2} + \\ &\qquad\qquad\qquad\qquad
         -B^+_{j+1,-m_1+1}\hat{C}^+_{l+1,j+1,-m_2}B^-_{j+1,m_1-1}\hat{C}^+_{l,j,m_2}\bigr\}\ket{l+2,m_1-2,m_2;j+2} \;,\\
&\oh[x_1^*-\delta(x_1^*),Jx_1J]\ket{l,m_1,m_2;j}  \\
&\qquad\simeq \bigl\{B^+_{j,m_1-2}\hat{C}^+_{l,j,m_2}B^+_{j,-m_1}\hat{C}^+_{l,j,-m_2} + \\ &\qquad\qquad\qquad\qquad
          -B^+_{j-1,-m_1+1}\hat{C}^+_{l-1,j-1,-m_2}B^+_{j-1,m_1-1}\hat{C}^+_{l-1,j-1,m_2}\bigr\}\ket{l,m_1-2,m_2;j} \\
&\qquad\qquad +\bigl\{B^+_{j-2,m_1-2}\hat{C}^+_{l-2,j-2,m_2}B^-_{j,-m_1}\hat{C}^+_{l-1,j-1,-m_2} + \\ &\qquad\qquad\qquad\qquad
          -B^-_{j-1,-m_1+1}\hat{C}^+_{l-2,j-2,-m_2}B^+_{j-1,m_1-1}\hat{C}^+_{l-1,j-1,m_2}\bigr\}\ket{l-2,m_1-2,m_2;j-2} \;.
\end{align*}
We need to estimate products of the form $\hat{C}^+_{l,j,-m_2}\hat{C}^+_{l+i,j+i,m_2}$,
for which, modulo smoothing operators, we find
\begin{align*}
\hat{C}^+_{l,j,-m_2}\hat{C}^+_{l+i,j+i,m_2}
&\simeq\sqrt{q^{2(l-j)}-q^{2l}q^{2(l-m_2+3-\epsilon)}}\sqrt{q^{2(l-j)}-q^{2l}q^{2(l+m_2+3+2i+\epsilon)}} \\
&\simeq \sqrt{q^{2(l-j)}}\sqrt{q^{2(l-j)}}=q^{2(l-j)} \;.
\end{align*}
Using this we get
\begin{align*}
&\oh[x_1+\delta(x_1),Jx_0J]\ket{l,m_1,m_2;j} \\
&\qquad \simeq q^{2(l-j)}\bigl\{A^+_{j+1,-m_1-1}B^+_{j,m_1}
         -A^+_{j,-m_1}B^+_{j+1,m_1}\bigr\}\ket{l+2,m_1+1,m_2;j+2} \\
&\qquad +q^{2(l-j)}\bigl\{A^+_{j,-m_1-1}B^+_{j,m_1}
         -A^+_{j-1,-m_1}B^+_{j-1,m_1}\bigr\}\ket{l,m_1+1,m_2;j} \;, \\
&\oh[x_1+\delta(x_1),Jx_1J]\ket{l,m_1,m_2;j} \\
&\qquad \simeq q^{2(l-j)}\bigl\{B^+_{j+1,m_1-1}B^+_{j,-m_1}
         -B^+_{j+1,-m_1-1}B^+_{j,m_1}\bigr\}\ket{l+2,m_1,m_2;j+2} \\
&\qquad +q^{2(l-j)}\bigl\{B^+_{j-1,m_1-1}B^-_{j,-m_1}
         -B^-_{j+1,-m_1-1}B^+_{j,m_1}\bigr\}\ket{l,m_1,m_2;j} \;,\\
&\oh[x_1-\delta(x_1),Jx_1J]\ket{l,m_1,m_2;j} \\
&\qquad \simeq q^{2(l-j)}\bigl\{B^-_{j+1,m_1-1}B^+_{j,-m_1}-B^-_{j,m_1}B^+_{j-1,-m_1-1}\bigr\}\ket{l,m_1,m_2;j} \;,\\
&\oh[x_1^*+\delta(x_1^*),Jx_1J]\ket{l,m_1,m_2;j} \\
&\qquad\simeq q^{2(l-j)}\bigl\{B^-_{j+2,m_1-2}B^+_{j,-m_1}-B^+_{j+1,-m_1+1}B^-_{j+1,m_1-1}\bigr\}\ket{l+2,m_1-2,m_2;j+2} \;,\\
&\oh[x_1^*-\delta(x_1^*),Jx_1J]\ket{l,m_1,m_2;j}  \\
&\qquad\simeq q^{2(l-j)}\bigl\{B^+_{j,m_1-2}B^+_{j,-m_1}-B^+_{j-1,-m_1+1}B^+_{j-1,m_1-1}\bigr\}\ket{l,m_1-2,m_2;j} \\
&\qquad +q^{2(l-j)}\bigl\{B^+_{j-2,m_1-2}B^-_{j,-m_1}-B^-_{j-1,-m_1+1}B^+_{j-1,m_1-1}\bigr\}\ket{l-2,m_1-2,m_2;j-2} \;.
\end{align*}
To prove that these commutators are smoothing we still need to check
that the terms in braces are bounded by $q^j$ (since $q^{2(l-j)}q^j\leq q^l$
is of rapid decay). This is done by using Eqs.~(\ref{eq:get}).
For example the first two braces are identically
zero, while the third one is
\begin{align*}
B^+_{j+1,m_1-1}B^+_{j,-m_1}-B^+_{j+1,-m_1-1}B^+_{j,m_1}
&=\tilde{B}^+_{j+1,m_1-1}\tilde{B}^+_{j,-m_1}-\tilde{B}^+_{j+1,-m_1-1}\tilde{B}^+_{j,m_1}
+O(q^j) \\ &=0+O(q^j) \;.
\end{align*}

What remains to control are the commutators $[x_2+\delta(x_2),JbJ]$
for $b=x_0,x_1,x_2$ and the commutators $[x_2^*-\delta(x_2^*),JbJ]$
for $b=x_1,x_2$ (which involve the `doubly-boxed' term).

The operators $x_2$ and $\delta(x_2)$ do not shift $m_1,j$ and have
coefficients independent on $m_1$. Thus, any operator acting only on
the label $m_1$ and with coefficients depending only on $m_1,j$,
commutes with $x_2$ and $\delta(x_2)$ and so can be neglected.
In particular, $x_0$ and $x_1$ can be written as sums of products of
operators of this kind (commuting with $x_2$ and $\delta(x_2)$) by
operators $y_i$'s,
\begin{align*}
y_1\ket{l,m_1,m_2;j} & :=\hat{C}^+_{l,j,m_2}\ket{l+1,m_1,m_2;j+1} \;,\\
y_2\ket{l,m_1,m_2;j} & :=\hat{H}^+_{l,j,m_2}\ket{l+1,m_1,m_2;j} \;,\\
y_3\ket{l,m_1,m_2;j} & :=\hat{C}^-_{l,j,m_2}\ket{l-1,m_1,m_2;j+1} \;,
\end{align*}
and their adjoints.
To prove that the commutators
$[x_2+\delta(x_2),JbJ]$,
for $b=x_0,x_1,x_2$, and $[x_2^*-\delta(x_2^*),JbJ]$,
for $b=x_1,x_2$, are smoothing, is sufficient to establish
the same for $b=y_1,y_2,y_3$.
For these operators we have
\begin{align*}
\oh[x_2+\delta(x_2),&\,Jy_1J]\ket{l,m_1,m_2;j} \\ & \simeq
        \bigl\{\hat{C}^+_{l+1,j,-m_2-1}\hat{D}_{l,j,m_2}^+
        -\hat{D}_{l+1,j+1,m_2}^+\hat{C}^+_{l,j,-m_2}\bigr\}
        \ket{l+2,m_1,m_2+1;j+1} \\ &=
        \hat{C}^+_{l,j,-m_2}\bigl\{\hat{D}_{l,j,m_2}^+
        -\hat{D}_{l+1,j+1,m_2}^+\bigr\}\ket{l+2,m_1,m_2+1;j+1} \;,\\
\oh[x_2+\delta(x_2),&\,Jy_2J]\ket{l,m_1,m_2;j} \\ & \simeq
        \bigl\{\hat{H}^+_{l+1,j,-m_2-1}\hat{D}_{l,j,m_2}^+
        -\hat{D}_{l+1,j,m_2}^+\hat{H}^+_{l,j,-m_2}\bigr\}
        \ket{l+2,m_1,m_2+1;j} \\ &=
        \hat{H}^+_{l,j,-m_2}\bigl\{\hat{D}_{l,j,m_2}^+
        -\hat{D}_{l+1,j,m_2}^+\bigr\}\ket{l+2,m_1,m_2+1;j} \;,\\
\oh[x_2+\delta(x_2),&\,Jy_3J]\ket{l,m_1,m_2;j} \\ & \simeq
        \bigl\{\hat{C}^-_{l+1,j,-m_2-1}\hat{D}_{l,j,m_2}^+
        -\hat{D}_{l-1,j+1,m_2}^+\hat{C}^-_{l,j,-m_2}\bigr\}
        \ket{l,m_1,m_2+1;j+1} \\ &=
        \hat{C}^-_{l,j,-m_2}\bigl\{\hat{D}_{l,j,m_2}^+
        -\hat{D}_{l-1,j+1,m_2}^+\bigr\}\ket{l,m_1,m_2+1;j+1} \;,\\
\oh[x_2^*-\delta(x_2^*),&\,Jy_1J]\ket{l,m_1,m_2;j} \\ & \simeq
        \bigl\{\hat{C}^+_{l-1,j,-m_2+1}\hat{D}_{l-1,j,m_2-1}^+
        -\hat{D}_{l,j+1,m_2-1}^+\hat{C}^+_{l,j,-m_2}\bigr\}
        \ket{l,m_1,m_2-1;j+1} \\ &=
        \hat{C}^+_{l,j,-m_2}\bigl\{\hat{D}_{l-1,j,m_2-1}^+
        -\hat{D}_{l,j+1,m_2-1}^+\bigr\}\ket{l,m_1,m_2-1;j+1} \;,\\
\oh[x_2^*-\delta(x_2^*),&\,Jy_2J]\ket{l,m_1,m_2;j} \\ & \simeq
        \bigl\{\hat{H}^+_{l-1,j,-m_2+1}\hat{D}_{l-1,j,m_2-1}^+
        -\hat{D}_{l,j,m_2-1}^+\hat{H}^+_{l,j,-m_2}\bigr\}
        \ket{l,m_1,m_2-1;j} \\ &=
        \hat{H}^+_{l,j,-m_2}\bigl\{\hat{D}_{l-1,j,m_2-1}^+
        -\hat{D}_{l,j,m_2-1}^+\bigr\}\ket{l,m_1,m_2-1;j} \;,\\
\oh[x_2^*-\delta(x_2^*),&\,Jy_3J]\ket{l,m_1,m_2;j} \\ & \simeq
        \bigl\{\hat{C}^-_{l-1,j,-m_2+1}\hat{D}_{l-1,j,m_2-1}^+
        -\hat{D}_{l-2,j+1,m_2-1}^+\hat{C}^-_{l,j,-m_2}\bigr\}
        \ket{l-2,m_1,m_2+1;j+1} \\ &=
        \hat{C}^-_{l,j,-m_2}\bigl\{\hat{D}_{l-1,j,m_2-1}^+
        -\hat{D}_{l-2,j+1,m_2-1}^+\bigr\}\ket{l-2,m_1,m_2+1;j+1} \;.
\end{align*}
Now $\hat{D}_{l\pm 1,j+1,m_2}^+-\hat{D}_{l,j,m_2}^+$ is bounded
by $q^{l\pm j+m_2}$, and $q^{l\pm j+m_2}\hat{C}^\pm_{l,j,-m_2}$
is bounded by $q^{2l}$. Next, $\hat{H}^+_{l,j,-m_2}$ is bounded by
$q^{l+j-m_2}$, and $q^{l+j-m_2}\hat{D}_{l,j,m_2}^+\simeq 1$.
This proves that all previous commutators are smoothing.
For the $y_i^*$'s the same statement follows from the symmetry (\ref{eq:symstar}).

We have arrived at last two commutators.
Modulo smoothing operators, the first one is
\begin{align*}
\oh[x_2+\delta(x_2),&\,Jx_2J]\ket{l,m_1,m_2;j} \\ & \simeq
\bigl\{\hat{D}_{l+1,j,-m_2-1}^+\hat{D}_{l,j,m_2}^+
-\hat{D}_{l+1,j,m_2-1}^+\hat{D}_{l,j,-m_2}^+\bigr\}\ket{l+2,m_1,m_2;j} \\
&\qquad\qquad +\bigl\{\hat{D}_{l+1,j,-m_2-1}^-\hat{D}_{l,j,m_2}^+
-\hat{D}_{l-1,j,m_2-1}^+\hat{D}_{l,j,-m_2}^-\bigr\}\ket{l,m_1,m_2;j} \\
&=\hat{D}_{l,j,-m_2}^-\bigl\{\hat{D}_{l,j,m_2}^+
-\hat{D}_{l-1,j,m_2-1}^+\bigr\}\ket{l,m_1,m_2;j} \;, \\
\intertext{where the second equality follows from the fact
that both $\hat{D}_{l,j,m_2}^+$ and $\hat{D}_{l,j,m_2}^-$ in
(\ref{eq:Dhat}) depend on $l$ and $m_2$ only through their sum.
For the same reason we have also that}
\oh[x_2^*-\delta(x_2^*),&\,Jx_2J]\ket{l,m_1,m_2;j} \\ & \simeq
\bigl\{\hat{D}_{l-1,j,-m_2+1}^+\hat{D}_{l-1,j,m_2-1}^+
-\hat{D}_{l,j,m_2-2}^+\hat{D}_{l,j,-m_2}^+\bigr\}\ket{l,m_1,m_2-2;j} \\
&\qquad\qquad +\bigl\{\hat{D}_{l-1,j,-m_2+1}^-\hat{D}_{l-1,j,m_2-1}^+
-\hat{D}_{l-2,j,m_2-2}^+\hat{D}_{l,j,-m_2}^-\bigr\}\ket{l-2,m_1,m_2-2;j} \\
&=\hat{D}_{l,j,-m_2}^-\bigl\{\hat{D}_{l-1,j,m_2-1}^+
-\hat{D}_{l-2,j,m_2-2}^+\bigr\}\ket{l-2,m_1,m_2-2;j} \;.
\end{align*}
The final observation that $\hat{D}_{l,j,-m_2}^-\hat{D}_{l-i,j,m_2-i}^+\simeq\hat{D}_{l,j,-m_2}^-$,
for $i=0,1,2$, gives that these commutators vanish modulo smoothing operators.
\end{prova}


\subsection*{Acknowledgements}
We are grateful to the referee whose remarks led to a much improved version of the paper.
This work was partially supported by the `Italian project Cofin06 -- Noncommutative geometry,
quantum groups and applications'.


\end{document}